\documentclass[a4paper,12pt]{amsart}
\usepackage[cp1251]{inputenc}
\usepackage[T2A]{fontenc}
\usepackage[russian,ukrainian,english]{babel}
\usepackage{amscd}
\usepackage{amssymb}
\usepackage{amsmath,amsthm,amsfonts}
\usepackage{hhline}

\usepackage{cite}

\newlength\deltaH
\newlength\deltaHcenter
\newlength\deltaV
\newlength\deltaVcenter
\newtheorem{theorem}{Theorem}[section]
\newtheorem{lem}{Lemma}[section]

\newtheorem{ozn}{Definition}[section]
\newtheorem{remk}{Remark}[section]
\newtheorem{nas}{Corollary}[section]
\newtheorem{prop}{Proposition}[section]

\newcommand{\rr}{\mathbb{R}}
\newcommand{\nn}{\mathbb{N}}
\renewcommand{\emptyset}{\varnothing}
\newcommand{\Cl}[1]{\overline{#1}}
\newcommand{\Int}[1]{\mathrm{Int}\,{#1}}
\newcommand{\Fr}[1]{\mathrm{Fr}\,{#1}}
\newcommand{\cD}{\mathcal{D}}
\newcommand{\cP}{\mathcal{P}}
\newcommand{\ff}{\mathfrak{f}}
\newcommand{\pr}{\mathop\mathrm{pr}\nolimits}
\newcommand{\restrict}[2]{#1\raisebox{-0.65ex}{$\left|\vphantom{{#1}_{#2}}\right.$}_{#2}}

\usepackage[dvips]{graphicx}

\begin{document}

\begin{center}
{\large\bf On the criteria of $\cD$-planarity of a tree.}
\end{center}
\vspace{3mm}
\begin{center}
{\small Polulyakh E.,\;Yurchuk I.}
\end{center}
\begin{center}
{\small Institute of Mathematics of Ukrainian national academy of sciences, Kyiv}
\end{center}
\vspace{5mm}

{\noindent \small\textbf{Abstract.} Let $T$ be a tree with a fixed subset of vertices $V^{\ast}$ such that there is a cyclic order $C$ on it and all terminal vertices are contained in this set. Let $D^{2} = \{ (x,y) \in \rr^{2} \,|\, x^{2} + y^{2} \leq 1 \}$
be a closed oriented 2--dimensional disk. The tree $T$ is called $\cD$-planar if there exists an embedding $\varphi : T \rightarrow \rr^{2}$ which satisfies the following conditions $\varphi(T) \subseteq D^{2}$, $\varphi(T) \cap \partial D^{2} = \varphi(V^{\ast})$ and if $\sharp V^{\ast} \geq 3$ then a cyclic order $\varphi(C)$ of $\varphi(V^{\ast})$ coincides with a cyclic order which is generated by the orientation of $\partial D^{2} \cong S^{1}$.

We obtain a necessary and sufficient condition for $T$ to be $\cD$-planar.  }
\medskip


\section*{Introduction.}

It is known that any tree is planar (i.e., can be embedded into a plane), more over, in~\cite{Ikebe} it was proved that any rooted tree with $n$ vertices can be embedded as a plane spanning tree on $n$ points of a plane, with the root being mapped onto an arbitrary specified point of them. There are similar results for rooted star forests~\cite{KK}. Many authors were interested in different types of embedding of trees into a disk with extra conditions. For example, in~\cite{NC1,NC2} author used a special type of embedding of tree into a disk for constructing a Morse function on it.

In this paper we will consider an embedding of a tree $T$ into $D^2$ (i.e., a closed oriented 2--dimensional disk) such that a fixed subset of its vertices $V^{\ast}$ which containes all terminal vertices maps to a boundary of $D^2$ in a special way and $T\setminus V^{\ast}$ maps into $Int D^2$. 

The existence of such embedding is a part of a solution of one topological problem: the realization of a finite graph as invariant of pseudoharmonic functions defined on a disk. Namely, a continuous function $f:D^2\rightarrow \mathbb{R}$ is called a pseudoharmonic function if $f\mid _{\partial D^2}$ is a continuous function with a finite number of local extrema and $f\mid_{\mathrm{Int} D^2}$ has a finite number of critical points and each of them is a saddle point (in the neighborhood of it $f$ has a representation as $Re z^n + const$, $z=x+iy$ and $n\geq 2$). In~\cite{Yu} the topological invariant of such functions is constructed, in particular, it consists of a cycle $\gamma$ with special properties and the complement to it is a disjoint union of trees. The criteria of topological equivalence of such functions is formulated in terms of their invariant.

In Section 1 we will describe some properties of a $\cD$-planar tree. They will be useful for the proof of Theorem 2.1 which gives the criteria of $\cD$-planarity of a tree. In this section we also study some properties of different types of relations on finite sets.

In Section 2 the criteria of $\cD$-planarity of a tree will be proved. Its proof has a topological nature.

We are sincerely grateful to V.V. Sharko, S.I. Maksymenko and I. Yu. Vlasenko for useful discussions. 



\section{The utility results.}

\subsection{Properties of trees embedded into two-dimensional disk.}

 Let $T$ be a tree with a set of vertices $V$ and a set of edges $E$. Suppose that $T$ is non degenerated ( has at least one edge). Denote by $V_{ter}$ a set of all vertices of $T$ such that their degree equals to 1. Let us assume that for a subset $V^{\ast} \subseteq V$ the following  condition holds true
\begin{equation}\label{eq_01}
V_{ter} \subseteq V^{\ast} \,.
\end{equation}
Let also $\varphi : T \rightarrow \rr^{2}$ is an embedding such that
\begin{equation}\label{eq_02}
\varphi(T) \subseteq D^{2} \,, \quad \varphi(T) \cap \partial D^{2} = \varphi(V^{\ast}) \,.
\end{equation}

\begin{lem}\label{lemma_01}
A set $\rr^{2} \setminus (\varphi(T) \cup \partial D^{2})$ has a finite number of connected components
\[
U_{0} = \rr^{2} \setminus D^{2}, U_{1}, \ldots, U_{m} \,,
\]
and for every $i \in \{1, \ldots, m\}$ a set $U_{i}$ is an open disk and is bounded by a simple closed curve
\[
\partial U_{i} = L_{i} \cup \varphi(P(v_{i}, v_{i}')) \,, \quad
L_{i} \cap \varphi(P(v_{i}, v_{i}')) = \{ \varphi(v_{i}), \varphi(v_{i}') \}
\]
 where $L_{i}$ is an arc of $\partial D^{2}$ such that the vertices $\varphi(v_{i})$ and $\varphi(v_{i}')$ are its endpoints, and  $\varphi(P(v_{i}, v_{i}'))$ is an image of the unique path $P(v_{i}, v_{i}')$ in $T$ which connects $v_{i}$ and $v_{i}'$.
\end{lem}

\begin{proof}
We prove lemma by an induction on the number of elements of the set $V^{\ast}$. Denote by $\sharp A$ a number of elements of a set $A$ .

Let us remark that $\sharp V^{\ast} \geq 2$ since $V_{ter} \subseteq V^{\ast}$ and $\sharp V_{ter} \geq 2$ for a non degenerated tree $T$ (it is easily verified by induction on the number of vertices).

\medskip
{\slshape Base of induction.} Let $\sharp V^{\ast} = 2$. From what was said above it follows that $V^{\ast} = V_{ter}$. So, a  tree satisfies a condition $\sharp V_{ter} = 2$. For such trees it is easy to prove by induction on the number of vertices that every vertex of $V \setminus V_{ter}$ has degree 2. In other words it is adjacent to two edges.

If a tree is considered as CW-complex (i.e. 0-cells are its vertices and 1-cells are its edges), then a topological space $T$ is homeomorphic to a segment with a set of the endpoints which coincides with $V^{\ast} = V_{ter}$.

Let $\varphi(T)$ be a cut of a disk $D^{2}$ between $\varphi(v_{1})$ and $\varphi(v_{2})$, where $\{ v_{1}, v_{2} \} = V_{ter}$.

Let us fix a homeomorphism
\[
\Phi_{0} : \partial D^{2} \cup \varphi(T) \rightarrow
\partial D^{2} \cup ( [-1, 1] \times \{0\} ) \,,
\]
such that $\Phi_{0} \circ \varphi (v_{1}) = (-1, 0)$, $\Phi_{0} \circ \varphi (v_{2}) = (1, 0)$, $\Phi_{0}(\partial D^{2}) = \partial D^{2}$, $\Phi_{0} \circ \varphi (T) = [-1, 1] \times \{ 0 \}$.

By Shernflic's theorem~\cite{Newman,Z-F-C} we can find a homeomorphism $\Phi : \rr^{2} \rightarrow \rr^{2}$ which extends $\Phi_{0}$. It is obvious that an embedding $\Phi \circ \varphi : T \rightarrow \rr^{2}$ complies with the conditions of lemma. From fact that $\Phi$ is homeomorphism it follows that $\varphi$ satisfies conditions of lemma.

\medskip
{\slshape Step of induction.}
Suppose that for some $n > 2$ lemma is proved for all trees with $\sharp V^{\ast} < n$ and their embeddings into $\rr^{2}$ which hold Conditions~\eqref{eq_01} and~\eqref{eq_02}.

Let a tree $T$ such that $V_{ter} \subseteq V^{\ast}$, $\sharp V^{\ast} = n$, and an embedding $\varphi : T \rightarrow \rr^{2}$ which satisfy Conditions~\eqref{eq_02} is fixed.

As we noticed above the set $V_{ter}$ contains at least two elements $w_{1}, w_{2} \in V_{ter}$. Let us consider the path $P(w_{1}, w_{2})$ which connects those vertices. Suppose that it passes through the vertices in the following order $w_{1} = u_{0}, u_{1}, \ldots, u_{k-1}, u_{k} = w_{2}$. Every vertex $u_{1}, \ldots, u_{k-1}$ has degree at least 2 since it is adjacent to two edges of $P(w_{1}, w_{2})$.

There exists a vertex $u_{s}$, $s \in \{1, \ldots, k-1\}$ such that
\begin{enumerate}
    \item[(i)] a degree of $u_{i}$ equals to 2 and $u_{i} \notin V^{\ast}$ for $i \in \{1, \ldots, s-1\}$;
    \item[(ii)] either a degree of $u_{s}$ is greater than 2 or $u_{s} \in V^{\ast}$ and a degree of $u_{s}$ equals to 2.
\end{enumerate}

Remark that a degree of $u_{s}$ does not equal to 1. Otherwise, the correlations $u_{s} = w_{2}$, $T = P(v_{1}, v_{2})$, $V^{\ast} = \{w_{1}, w_{2}\}$, $\sharp V^{\ast} = 2$ should be satisfied but we assumed that $\sharp V^{\ast} \geq 3$.

Let us consider a path $P(w_{1}, u_{s}) = P(u_{0}, u_{s})$. Suppose that it passes through edges $e_{1}, \ldots, e_{s}$ successively.

We consider a subgraph $T'$ of $T$ with the set of vertices and edges, respectively, as followes
\[
V(T') = V \setminus \{u_{0}, \ldots, u_{s-1}\} \,, \quad
E(T') = E \setminus \{e_{1}, \ldots, e_{s}\} \,.
\]
By construction $u_{0} \in V_{ter}(T)$ and $u_{0}$ is adjacent to $e_{1}$ in $T$; every vertex $u_{i}$, $i \in \{1, \ldots , s-1\}$ has degree 2 thus it is adjacent only to $e_{i}$ and $e_{i+1}$ in $T$. Therefore a graph $T'$ is defined correctly.

A graph $T'$ has no cycles since it is a subgraph of $T$. Let us verify that $T'$ is connected. Let $v', v'' \in V(T')$ and  $P(v', v'')$ be a path which connects vertices $v'$ and $v''$ in $T$. Then a path $P(v', v'')$ does not pass through a vertex $u_{0} = w_{1}$ since $u_{0} \in V_{ter}$ and only one edge $e_{1}$ is adjacent to this vertex. Thus $e_{1} \notin P(v', v'')$. Similarly, if $s \geq 2$ then $e_{2} \notin P(v', v'')$ since  an edge $e_{2}$ is adjacent to a vertex $u_{1}$ which is in addition adjacent only to $e_{1}$ and $e_{1} \notin P(v', v'')$. Similarly, by induction we prove that $e_{i} \notin P(v', v'')$ for every $i \in \{1, \ldots, s\}$. Thus a path $P(v', v'')$ connects vertices $v'$ and $v''$ in $T'$. Therefore a graph  $T'$ is connected.

We verified that $T'$ is a tree. Let us define $V^{\ast}(T') = V^{\ast}(T) \cap V(T')$, $\varphi_{0} = \varphi |_{T'} : T' \rightarrow \rr^{2}$. By definition of a set $V^{\ast}(T')$ it is obvious that a map $\varphi_{0}$ satisfies condition~\eqref{eq_02}. Also $\sharp V^{\ast}(T') < \sharp V^{\ast}(T)$ since $u_{0} \in V^{\ast}(T) \setminus V^{\ast}(T')$. Thus $\sharp V^{\ast}(T') < n$.

Let us check that $V_{ter}(T') \subseteq V^{\ast}(T')$.

By construction for every vertex $v \neq u_{s}$ of $T'$ its degrees coincide in $T$ and $T'$. The degree of $u_{s}$ in $T'$ is on one less then degree of $u_{s}$ in $T$. Thus $V_{ter}(T') \subseteq V_{ter}(T) \cup \{u_{s}\}$.

If $u_{s} \in V^{\ast}(T)$, then $V_{ter}(T') \subseteq V_{ter}(T) \cup V^{\ast}(T) \subseteq V^{\ast}(T)$. Therefore $V_{ter}(T') \subseteq V^{\ast}(T) \cap V(T') = V^{\ast}(T')$.

Let $u_{s} \notin V^{\ast}(T)$. By definition the degree of $u_{s}$ in $T$ is not less then 3 and a degree of $u_{s}$ in $T'$ is not less then 2. Thus $V_{ter}(T') \subseteq V_{ter}(T) \subseteq V^{\ast}(T)$. So, as above, $V_{ter}(T') \subseteq V^{\ast}(T')$.

By induction lemma holds true for a tree $T'$ and an embedding $\varphi_{0} : T' \rightarrow \rr^{2}$.

Denote by $W_{0} = \rr^{2} \setminus D^{2}, W_{1}, \ldots, W_{r}$ connected components of a set $\rr^{2} \setminus (\varphi_{0}(T') \cup \partial D^{2})$.

It is obvious that
\[
\varphi(T) = \varphi(T') \cup \varphi(P(u_{0}, u_{s})) =
\varphi_{0}(T') \cup \varphi(P(u_{0}, u_{s})) \,.
\]
Therefore $\varphi(T) \cup \partial D^{2} = (\varphi_{0}(T') \cup \partial D^{2}) \cup \varphi(P(u_{0}, u_{s}))$.
By construction we get that$(\varphi(T') \cup \partial D^{2}) \cap \varphi(P(u_{0}, u_{s})) = \{\varphi(u_{0}), \varphi(u_{s})\}$.

Denote $J = \varphi(P(u_{0}, u_{s}))$. The set $J_{0} = J \setminus \{\varphi(u_{0}), \varphi(u_{s})\}$ is a homeomorphic image of interval thus it is connected. But besides  $J_{0} \cap (\varphi_{0}(T') \cup \partial D^{2}) = \emptyset$ thus there exists a component $W_{j}$ which contains  $J_{0}$ (it is easy to see that $j \neq 0$).

By assumption of induction the boundary of disk $W_{j}$ is a simple closed curve $\partial W_{j} = K_{j} \cup \varphi_{0}(P(v_{j}, v_{j}'))$ which consists of an arc $K_{j}$ of a circle $\partial D^{2}$ with the ends $\varphi_{0}(v_{j})$ and $\varphi_{0}(v_{j}')$ and an image of path $P(v_{j}, v_{j}')$ which connects vertices $v_{j}, v_{j}' \in V^{\ast}(T')$ in $T'$ (this path also connects vertices $v_{j}$ and $v_{j}'$ in $T$).

 The set $J$ is a homeomorphic image of segment and also $J_{0} \subseteq W_{j}$, $\varphi(u_{0}) \in \partial D^{2} \subseteq (\rr^{2} \setminus W_{j})$, $\varphi(u_{s}) \in \varphi_{0}(T') \subseteq (\rr^{2} \setminus W_{j})$. Therefore $J$ is a cut of disk $W_{j}$ between points $\varphi(v_{j})$ and $\varphi(v_{j}')$. Correlations $\varphi(u_{s}) \in \varphi(P(v_{j}, v_{j}'))$, $\varphi(u_{0}) \in K_{j} \setminus \{\varphi(v_{j}), \varphi(v_{j}')\} = \partial W_{j} \setminus \varphi(T')$ hold true since $u_{0} \notin V(T')$ and $\varphi(u_{0}) \notin \varphi(T')$.

So, a set $\Cl{W}_{j} \setminus ( \partial W_{j} \cup \varphi(P(u_{0}, u_{s})) )$ has two connected components $W_{j}^{1}$, $W_{j}^{2}$ which are homeomorphic to open disks and bounded by simple closed curves.

We remark that the arc $\varphi(P(v_{j}, v_{j}'))$ is not a point, otherwise the correlations $K_{j} \cong \partial D^{2}$, $\varphi_{0}(T') \cap \partial D^{2} = \{ \varphi(v_{j}) = \varphi(v_{j}') \}$, $\sharp V^{\ast}(T') = \sharp (\varphi_{0}(T') \cap \partial D^{2}) = 1$ should hold true. Thus points $\varphi(v_{j})$ and $\varphi(v_{j}')$ are different. From the inclusions  $\varphi(u_{s}) \in \varphi(P(v_{j},v_{j}'))$, $\varphi(u_{0}) \in \partial W_{j} \setminus \varphi(P(v_{j}, v_{j}'))$ it follows that points $\varphi(v_{j})$ and $\varphi(v_{j}')$ can not be contained in a set $\partial W_{j}^{1} \cap \partial W_{j}^{2} = \varphi(P(u_{0},u_{s}))$ simultaneously.

Let $\varphi(v_{j}) \in \partial W_{j}^{1}$, $\varphi(v_{j}') \in \partial W_{j}^{2}$. By those correlations the sets  $W_{j}^{1}$ and $W_{j}^{2}$ are defined uniquely.

Points $\varphi(u_{0})$, $\varphi(u_{s})$ divide the circle onto two arcs $R_{1}$, $R_{2}$ with $R_{1} \subseteq \partial W_{j}^{1} \setminus W_{j}^{2}$, $R_{2} \subseteq \partial W_{j}^{2} \setminus W_{j}^{1}$.

Suppose for some edge $e \in E(T)$ its image is contained in  $\partial W_{j}$. Then the image of $e$ without the ends is connected set and belongs to $\partial W_{j} \setminus \{\varphi(u_{0}), \varphi(u_{s})\} = R_{1} \cup R_{2}$. Thus the image of $e$ without the endpoints belongs to either $R_{1}$ or $R_{2}$.

The path  which connects vertices $v_{j}$ and $v_{j}'$ in $T'$ passes through the vertices $v_{j} = \hat{v}_{0}, \hat{v}_{1}, \ldots, \hat{v}_{k} = v_{j}'$ and through the edges $\hat{e}_{1}, \ldots, \hat{e}_{k}$ in this order.

If $\varphi(\hat{v}_{i}) \in R_{1}$ for some $i \in \{0, \ldots, k\}$, then $\varphi(\hat{v}_{i}) \in \rr^{2} \setminus \Cl{R}_{2}$ and $(\varphi(\hat{e}_{i}) \setminus \{\varphi(\hat{v}_{i}), \varphi(\hat{v}_{i+1})\}) \cap (\rr^{2} \setminus \Cl{R}_{2}) \neq \emptyset$ since a point  $\varphi(\hat{v}_{i})$ is a boundary for the set $\varphi(\hat{e}_{i}) \setminus \{\varphi(\hat{v}_{i}), \varphi(\hat{v}_{i+1})\}$ but $\rr^{2} \setminus \Cl{R}_{2}$ is an open neighborhood of this point. From what we said it follows that $\varphi(\hat{e}_{i}) \setminus \varphi(\hat{v}_{i+1}) \subseteq R_{1}$. Therefore $\varphi(\hat{v}_{i+1}) \in \Cl{R}_{1} = R_{1} \cup \{\varphi(u_{0}), \varphi(u_{s})\}$. Indeed, either $\varphi(\hat{v}_{i+1}) \in R_{1}$ or $\varphi(\hat{v}_{i+1}) = \varphi(u_{s})$ (and $\hat{v}_{i+1} = u_{s}$) since $u_{0} \notin V(T')$ by construction.

By assumption of induction $\varphi(u_{s}) \in \varphi(P(v_{j}, v_{j}')) = \varphi(\hat{v}_{0}, \hat{v}_{k}))$. Therefore $u_{s} \in \{\hat{v}_{0}, \ldots, \hat{v}_{k}\}$ and there exists an index $k_{0} \in \{0, \ldots, k\}$ such that $u_{s} = \hat{v}_{k_{0}}$.

The inductive application of our previous argument leads us to correlations $\varphi(P(\hat{v}_{0}, u_{s})) \setminus \varphi(u_{s}) = \varphi(P(v_{j}, u_{s})) \setminus \varphi(u_{s}) \subseteq R_{1}$ (in the case when $v_{j} = u_{s}$ we get $\varphi(P(v_{j}, u_{s})) = \varphi(u_{s})$).

Similar argument give $\varphi(P(u_{s}, v_{j}')) \setminus \varphi(u_{s}) \subseteq R_{2}$.

Finally we get
$\partial W_{j}^{1} = R_{1} \cup \varphi(P(u_{0}, u_{s})) = R_{1} \cup J$, $\partial W_{j}^{2} = R_{2} \cup J$; $\partial W_{j}^{1} \cap \varphi(P(v_{j}, v_{j}')) = \partial W_{j}^{1} \cap (\varphi(P(v_{j}, u_{s})) \cup \varphi(P(u_{s}, v_{j}'))) = \varphi(P(v_{j}, u_{s}))$;
$\partial W_{j}^{2} \cap \varphi(P(v_{j}, v_{j}')) = \varphi(P(u_{s}, v_{j}'))$.

Therefore $\varphi(T) \cap \partial W_{j}^{1} = (\varphi(T') \cup J) \cap \partial W_{j}^{1} = (\varphi(P(v_{j}, v_{j}')) \cup J) \cap \partial W_{j}^{1} = \varphi(P(v_{j}, u_{s})) \cup \varphi(P(u_{0}, u_{s})) = \varphi(P(v_{j}, u_{0}))$; $\varphi(T) \cap \partial W_{j}^{2} = \varphi(P(v_{j}', u_{0}))$.

It is easy to see that $\varphi(v_{j}) \neq \varphi(u_{0})$ and $\varphi(v_{j}') \neq \varphi(u_{0})$ since $v_{j}$, $v_{j}' \in V(T')$ but $u_{0} \notin V(T')$. Hence a set $\varphi(P(v_{j}, u_{0})) \setminus \{\varphi(v_{j}), \varphi(u_{0})\}$ is one of two connected components of the set $\partial W_{j}^{1} \setminus \{\varphi(v_{j}), \varphi(u_{0})\}$. Another connected component of this set is contained in $\partial W_{j} \setminus \varphi(T') = K_{j} \subseteq \partial D^{2}$ thus it is an arc of circle $\partial D^{2}$ which connects points $\varphi(v_{j})$ and $\varphi(u_{0})$. Denote it by $K_{j}^{1}$.

Similarly, $\partial W_{j}^{2} = \varphi(P(v_{j}', u_{0})) \cup K_{j}^{2}$, where $K_{j}^{2}$ is an arc of $\partial D^{2}$ which connects points $\varphi(v_{j}')$ and $\varphi(u_{0})$.

We proved that a compliment $\rr^{2} \setminus (\varphi(T) \cup \partial D^{2})$ has a finite number of connected components
\[
\rr^{2} \setminus D^{2} = W_{0}, W_{1}, \ldots, W_{j-1},
W_{j}^{1}, W_{j}^{2}, W_{j+1}, \ldots, W_{r} \,;
\]
and the components $W_{j}^{1}$ and $W_{j}^{2}$ satisfy the conditions of lemma. Finally we remark that the correlations $\partial W_{k} \cap \varphi(T) = \partial W_{k} \cap \varphi(T') = \partial W_{k} \cap \varphi_{0}(T')$ hold true for  $k > 0$, $k \neq j$ thus
\[
\partial W_{k} = K_{k} \cup \varphi_{0}(P(v_{k}, v_{k}')) = K_{k} \cup \varphi(P(v_{k}, v_{k}'))
\]
and the component $W_{k}$ satisfies lemma.

\end{proof}

\begin{nas}\label{nas_01}
Let $T$ be a tree with fixed subset of vertices $V^{\ast} \supseteq V_{ter}$ and $\varphi : T \rightarrow \rr^{2}$ an embedding which satisfies~\eqref{eq_02}.

Then the following conditions hold true.

1)In notation of Lemma~\ref{lemma_01}
\[
L_{i} \cap \varphi(T) = \{\varphi(v_{i}), \varphi(v_{i}')\} \,,
\quad i = 1, \ldots, m \,.
\]

2) If there exists an arc $L$ of circle $\partial D^{2}$ with the ends $\varphi(u_{1})$, $\varphi(u_{2})$ such that $L \cap \varphi(T) = \{\varphi(u_{1}), \varphi(u_{2})\}$ for some $u_{1}$, $u_{2} \in V^{\ast}$, then  there exists $k \in \{1, \ldots, m\}$ such that $L \cup \varphi(P(u_{1}, u_{2})) = \partial U_{k}$ (then $L = L_{k}$, $u_{1} = v_{k}$, $u_{2} = v_{k}'$).
\end{nas}

\begin{proof}
1) Suppose that an arc $L_{i} \setminus \{\varphi(v_{i}), \varphi(v_{i}')\}$ contains a point $\varphi(v) \in  \varphi(T)$ for some $i \in \{1, \ldots, m\}$. Thus $v \in V^{\ast}$. Let  $e \in E(T)$ be an edge of graph $T$ which is adjacent to a vertex $v$ and $v' \in V$ be another end of the edge $e$.

A set $J_{0} = \varphi(e) \setminus \{\varphi(v), \varphi(v')\}$ is connected, $\varphi(v)$ is a boundary point of it, $W = \rr^{2} \setminus \varphi(P(v_{i}, v_{i}'))$ is an open neighborhood of a point $\varphi(v)$. Thus $J_{0} \cap W \neq \emptyset$ and $e \notin P(v_{i}, v_{i}')$. Hence $J_{0} \cap \varphi(P(v_{i}, v_{i}')) = \emptyset$. By the conditions of lemma also $J_{0} \cap \partial D^{2} = \emptyset$. A set $\varphi(P(v_{i}, v_{i}'))$ is a cut of closed disk $D^{2}$. Obviously, by construction a set $Q = \Cl{U}_{i} \setminus \varphi(P(v_{i}, v_{i}')) = U_{i} \cup (L_{i} \setminus \{\varphi(v_{i}), \varphi(v_{i}')\})$ is a connected component of the compliment $D^{2} \setminus \varphi(P(v_{i}, v_{i}'))$ which contains a point $\varphi(v)$. That point is a boundary point of the connected subset $J_{0}$ of a space $D^{2} \setminus \varphi(P(v_{i}, v_{i}'))$ therefore $J_{0} \subseteq Q$.

But $U_{i} \subseteq \rr^{2} \setminus \varphi(T)$, $L_{i} \subseteq \partial D^{2}$ and $J_{0} \subseteq \varphi(T) \setminus \varphi(V) \subseteq \varphi(T) \setminus \partial D^{2}$. Thus $J_{0} \cap Q \subseteq (J_{0} \cap U_{i}) \cup (J_{0} \cap L_{i}) = \emptyset$.

The contradiction obtained is a last step of the proof of first condition of corollary.

\medskip

2) Support that there exists an arc $L$ of $\partial D^{2}$ with the ends in points $\varphi(u_{1})$, $\varphi(u_{2})$ such that $L \cap \varphi(T) = \{\varphi(u_{1}), \varphi(u_{2})\}$ for some $u_{1}$, $u_{2} \in V^{\ast}$.

An arc $L$ bounders to some connected component $U_{k}$, $k \geq 1$ of the compliment $\rr^{2} \setminus (\varphi(T) \cup \partial D^{2})$. From Lemma~\ref{lemma_01} and first condition of corollary it follows that $\{u_{1}, u_{2}\} = \{v_{k}, v_{k}'\}$. Thus vertices $u_{1}$ and $u_{2}$ can be connected by a path $\tilde P(u_{1}, u_{2}) = P(v_{k}, v_{k'})$ which satisfies Lemma~\ref{lemma_01}. A graph $T$ is a tree thus  $P(u_{1}, u_{2}) = \tilde P(u_{1}, u_{2}) = P(v_{k}, v_{k'})$.
\end{proof}

Let $T$ be a tree with a fixed subset of vertices $V^{\ast}$ and $\varphi : T \rightarrow \rr^{2}$ is an embedding which satisfy~\eqref{eq_01} and~\eqref{eq_02}.

The pair of vertices $v_{1}$, $v_{2} \in V^{\ast}$, $v_{1} \neq v_{2}$ is said to be \emph{adjacent on a circle $\partial D^{2}$} if there exists an arc $L$ of this circle with the ends $\varphi(v_{1})$ and $\varphi(v_{2})$ such that $L \cap \varphi(T) = \{\varphi(v_{1}), \varphi(v_{2})\}$ holds true for it.

Denote by $\mathcal{P}$ a set of all paths in $T$ which connect adjacent pairs of vertices.

\begin{nas}\label{nas_02}
If $\sharp V^{\ast} \geq 3$, then a correspondence
\begin{gather*}
\Theta : \{U_{1}, \ldots, U_{m}\} \rightarrow \mathcal{P} \,,\\
\Theta(U_{i}) = P(v_{i}, v_{i}') \,,
\end{gather*}
is a bijective map.
\end{nas}

\begin{proof}
It is sufficient to check an injectivity of the map $\Theta$.

 Suppose that the following equalities hold true $\partial U_{i} = L_{i} \cup \varphi(P(v, v'))$, $\partial U_{j} = L_{j} \cup \varphi(P(v, v'))$  for some $i$, $j \in \{1, \ldots, m\}$, $i \neq j$. Then $L_{i} \cap L_{j} = \{\varphi(v), \varphi(v')\}$, $L_{i} \cup L_{j} \cong S^{1}$ therefore $L_{i} \cup L_{j} = \partial D^{2}$.

But from Corollary~\ref{nas_01} it follows that $\emptyset = (L_{i} \cup L_{j} \setminus \{\varphi(v), \varphi(v')\}) \cap \varphi(T)$. Therefore $\sharp (\partial D^{2} \cap \varphi(T)) = \sharp V^{\ast} \leq 2$ and it contradicts to the conditions of corollary.
\end{proof}


\subsection{On relations defined on finite sets.}

At first we remind that a ternary relation $O$ on the set $A$ is any subset of the $3^{rd}$ cartesian power $A^3:O\subseteq A^3$.

Let $A$ be a set, $O$ a ternary relation on $A$ which is asymmetric ($(x,y,z)\in O \Rightarrow (z,y,x)\overline{\in} O $), transitive $(x,y,z)\in O, (x,z,u)\in O \Rightarrow (x,y,u\in O)$ and cyclic $(x,y,z)\in O \Rightarrow (y,z,x)\in O$. Then $O$ is called a cyclic order on the set $A$~\cite{Nov}.

A cyclic order $O$ is a complete on a finite set $A$, $\sharp A \geq 3$, if $x,y,z\in A, x\neq y \neq z \neq x \Rightarrow $ there exists a permutation $(u,v,w)$ of sequence $(x,y, z)$ such that $(u,v,w) \in O$.

\begin{prop}\label{prop_01}
Let there is a complete cyclic order $O$ on some finite set $A$, $\sharp A \geq 3$.

Then for every $a \in A$ there exist unique $a'$, $a'' \in A$ such that
\begin{itemize}
    \item $O(a', a, b)$ for all $b \in A \setminus \{a, a'\}$;
    \item $O(a, a'', b)$ for all $b \in A \setminus \{a, a''\}$,
\end{itemize}
and $a' \neq a''$.
\end{prop}

\begin{proof}
Let us fix $a \in A$. By using~\cite{Nov} we can construct a binary relation $\rho$ up to the relation $O$ with the help of the following condition
\[
O(a, a_{1}, a_{2}) \Leftrightarrow a_{1} \,\rho\, a_{2} \,.
\]
It is easy to verify that the relation $\rho$ defines a strict linear order on a set $A \setminus \{a\}$.

 The set $A \setminus \{a\}$ is finite therefore there exist a minimal element $a'$ and maximal element $a''$ with respect to the order $\rho$ on this set. It is obvious that they satisfy conditions of proposition by definition.

Finally, $a' \neq a''$ since $\sharp (A \setminus \{a\}) \geq 2$.
\end{proof}

\begin{ozn}\label{ozn_cyclic_susidni}
Let there is a complete cyclic order $O$ on a set $A$, $\sharp A \geq 3$. Elements $a_{1}$, $a_{2} \in A$ are said to be \emph{adjacent} with respect to a cyclic order $O$ if one of the following conditions holds:
\begin{itemize}
    \item $O(a_{1}, a_{2}, b)$ for all $b \in A \setminus \{a_{1}, a_{2}\}$;
    \item $O(a_{2}, a_{1}, b)$ for all $b \in A \setminus \{a_{1}, a_{2}\}$.
\end{itemize}
\end{ozn}

\begin{remk}
From Proposition~\ref{prop_01} it follows that every element has exactly two adjacent elements on a finite set $A$ with a complete cyclic order.
\end{remk}


\begin{ozn}\label{ozn_zruchne}
Let $A$ be a finite set. A binary relation $\rho$ on $A$ is said to be \emph{convenient} if
\begin{itemize}
    \item[1)] for all $a$, $b \in A$ from $a \rho b$ it follows that $a \neq b$;
    \item[2)] for every $a \in A$ there is no more than one $a' \in A$ such that $a \rho a'$;
    \item[3)] for every $a \in A$ there is no more than one $a'' \in A$ such that $a'' \rho a$.
\end{itemize}
\end{ozn}

We remind that a \emph{graph} of the relation $\rho$ on $A$ is a set $\{(a,b) \in A \times A \,|\, a \rho b\}$.

Let $\rho$ be a convenient relation on a finite set $A$, $\hat\rho$ be a minimal relation of equivalence which contains $\rho$. Let us remind that a graph $\hat\rho$ consists of
\begin{itemize}
    \item  all pairs $(a,b)$ such that there exist $k = k(a,b) \in \nn$ and a sequence $a = a_{0}$, $a_{1}, \ldots, a_{k} = b$ which comply with one of the following conditions $a_{i-1} \rho a_{i}$, $a_{i} \rho a_{i-1}$ for every $i \in \{1, \ldots, k\}$;
    \item pairs $(a,a)$, $a \in A$.
\end{itemize}

We distinguished a diagonal $\Delta_{A \times A}$ since, in general, there could exist $a \in A$ such that neither $a \rho b$ nor $b \rho a$ holds true for all $b \in A$.

The relation $\hat\rho$ generates a partition $\ff$ of $A$ onto classes of equivalence.

\begin{prop}\label{prop_02}
Let $B \in \ff$ be a class of equivalence of the relation $\hat\rho$. Then there exists no more than one element $b \in B$ which is in the relation $\rho$ with no element of $A$.
\end{prop}

\begin{proof}
We remark that if either $a \rho b$ or $b \rho a$ and $b \in B$, then $a \in B$ by definition of $B$.

It is obvious that if $\sharp B = 1$ then proposition holds true. Let $\sharp B \geq 2$.

Let $a_{0}$, $a_{1}, \ldots, a_{k}$ be a fixed sequence of pairwise different elements of $B$ such that the correlation $a_{i-1} \rho a_{i}$ holds true  for any $i \in \{1, \ldots, k\}$.

If there exists $b \in B \setminus \{a_{0}, \ldots, a_{k}\}$, then there exists $b' \in B \setminus \{a_{0}, \ldots, a_{k}\}$ such that either $b' \rho a_{0}$ or $a_{k} \rho b'$. Let us verify it.

By definition of a set $B$ there exists a sequence $b = c_{0}$, $c_{1}, \ldots, c_{m} = a_{0}$ such that either  $c_{j-1} \rho c_{j}$ or $c_{j} \rho c_{j-1}$ holds true for all $j \in \{1, \ldots, m\}$. From correlations $c_{0} \notin \{a_{0}, \ldots, a_{k}\}$ and $c_{m} \in \{a_{0}, \ldots, a_{k}\}$ it follows that there is $s \in \{0, \ldots, m\}$ such that $c_{s-1} \notin \{a_{0}, \ldots, a_{k}\}$ but $c_{s} \in \{a_{0}, \ldots, a_{k}\}$. Thus $c_{s} = a_{r}$ for some $r \in \{0, \ldots, k\}$.

Let $c_{s-1} \rho c_{s}$, i.e. $c_{s-1} \rho a_{r}$. Then $r = 0$. Really, if $r \geq 1$, then $a_{r-1} \rho a_{r}$. By construction $c_{r-1} \neq a_{r-1}$ therefore a correlation $c_{s-1} \rho a_{r}$ contradicts to condition~3) of Definition~\ref{ozn_zruchne}.

Similarly, if $c_{s} \rho c_{s-1}$, then $c_{s} = a_{k}$.

It is easy to see that element $b' = c_{s-1}$ satisfies conditions of proposition.

From what we said above it follows that if for some pairwise different $a_{0}, \ldots, a_{k} \in B$ inequality $\{a_{0},  \ldots, a_{k}\} \neq B$ and relation $a_{i-1} \rho a_{i}$, $i \in \{1, \ldots, k\}$ hold true, then there are pairwise different $a_{0}', \ldots, a_{k+1}' \in B$ such that $a_{i-1}' \rho a_{i}'$, $i \in \{1, \ldots, k+1\}$ hold true for them.

By definition the set $B$ contains two elements $b'$, $b'' \in B$ such that $b' \rho b''$. So, by a finite number of steps (the set $B$ is finite) we can index all elements of $B$  in such way that the following correlations hold true
\begin{gather}\label{eq_max_chain}
a_{i-1} \rho a_{i} \,, \quad i \in \{1, \ldots, n\} \,; \\
\{a_{0}, \ldots, a_{n}\} = B \,.\nonumber
\end{gather}
Therefore only element $a_{n} \in B$ can satisfy conditions of the proposition.
\end{proof}

Let $\mu$ be some relation on a set $A$.

\begin{ozn}\label{ozn_cycle}
Elements $b_{0}, \ldots, b_{n} \in A$, $n \geq 1$ are said to generate \emph{$\mu$-cycle} if a graph of the relation $\mu$ contains a set
\begin{equation}\label{eq_cycle}
\{(b_{0}, b_{1}), \ldots, (b_{n-1}, b_{n}), (b_{n}, b_{0})\} \,.
\end{equation}
\end{ozn}

\begin{ozn}\label{ozn_chain}
Elements $b_{0}, \ldots, b_{n} \in A$, $n \geq 0$ to generate \emph{$\mu$-chain} if for arbitrary $a \in A$ the pairs $(a, b_{0})$ and $(b_{n}, a)$ do not belong to a graph of $\mu$ and for $n \geq 1$ a graph of the relation $\mu$ contains a set
\begin{equation}\label{eq_chain}
\{(b_{0}, b_{1}), \ldots, (b_{n-1}, b_{n})\} \,.
\end{equation}
\end{ozn}

\begin{nas}\label{nas_struct_1}
Let $\rho$ be a convenient relation, $B \in \ff$ a class of equivalence of the relation $\hat\rho$. Then the elements of $B$ generate either $\rho$-cycle or $\rho$-chain. In the first case a graph of the restriction of $\rho$ on the set $B$ is of form~\eqref{eq_cycle} and in the other it has form ~\eqref{eq_chain}.
\end{nas}

\begin{proof}
Let us order the elements of $B$ in such way that ~\eqref{eq_max_chain} holds true for them.

If there exists $a \in A$ such that $a_{n} \rho a$, then $a \in B$.
Conditions 1) and 3) of Definition~\ref{ozn_zruchne} obstruct to hold correlation $a_{i} \rho a_{j}$ for $i \neq j-1$, $j \in \{1, \ldots, n\}$.
Therefore $a = a_{0}$ and $a_{n} \rho a_{0}$.

Similarly, if there exists $a \in A$ which $a \rho a_{0}$, then from conditions 1) and 2) of Definition~\ref{ozn_zruchne} it follows that $a = a_{n}$ and $a_{n} \rho a_{0}$.

So, either a correlation $a_{n} \rho a_{0}$ holds true or for every $a \in A$ neither $a \rho a_{0}$ nor $a_{n} \rho a$ holds true. In the first case the elements of $B$ generate $\rho$-cycle (if $a_{n} \rho a_{0}$, then $a_{n} \neq a_{0}$ and $\sharp B \geq 2$ by definition), in the other case we get $\rho$-chain.
\end{proof}

\begin{nas}\label{nas_struct_2}
Let the elements of $B \subseteq A$ generate either $\rho$-cycle or $\rho$-chain. Then $B$ is a class of equivalence of the relation $\hat\rho$. If the elements of $B \subseteq A$ generate $\rho$-chain, then the relation $\rho$ generates a full linear order on $B$.
\end{nas}

\begin{proof}
Let $\hat\rho$ be a minimal relation of equivalence which contains $\rho$. By definition the set $B$ belongs to the unique class of equivalence  of the relation $\hat\rho$. Denote it by $\hat B$.

By definition the set $B$ satisfies~\eqref{eq_max_chain}. If there exists $b \in \hat B \setminus B$, then, as we verified in the proof of Proposition~\ref{prop_02}, there is $b' \in \hat B \setminus B$ such that
\begin{equation}\label{eq_06}
b' \rho a_{0} \quad \text{or} \quad a_{n} \rho b' \,.
\end{equation}
This contradicts to definition of $\rho$-chain. If the elements of $B$ generate $\rho$-cycle, then it follows from the definition of  convenient relation that

\begin{equation}\label{eq_07}
a_{n} \rho a_{0} \,,
\end{equation}
see Corollary~\ref{nas_struct_1}.
By using conditions 2) and 3) of a convenient relation  from equality  $b' \notin B$ we can conclude that~\eqref{eq_06} and~\eqref{eq_07} can not be satisfied simultaneously.

So, a set $B$ is a class of equivalence of the relation $\hat\rho$.

If elements of the set $B$ generate a chain, then a graph of a restriction of the relation $\rho$ on $B$ has form~\eqref{eq_chain}, see Corollary~\ref{nas_struct_1}. Therefore $\rho$ generates a linear order on the set $B$.
\end{proof}

\begin{ozn}\label{ozn_induc_by_cycle}
Let $O$ be a complete cyclic order on $A$, $\sharp A \geq 3$. $O$ is said to induce a binary relation $\rho_{O}$ on a $A$ according to the following rule:
$a \rho_{O} b$ if $O(a,b,c)$ $\forall c \in A \setminus \{a, b\}$.
\end{ozn}

From Proposition~\ref{prop_01} it follows that a relation $\rho_{O}$ is convenient.

\begin{prop}\label{prop_cycle}
If $O$ is a complete cyclic order on $A$, then all elements of $A$ generate $\rho_{O}$-cycle.
\end{prop}

\begin{proof}
Let $\hat\rho_{O}$ be a minimal relation of equivalence which contains $\rho_{O}$. From Proposition~\ref{prop_01} and Corollaries~\ref{nas_struct_1} and~\ref{nas_struct_2} it follows that every class of equivalence of the relation $\hat\rho_{O}$ is $\rho_{O}$-cycle and there are no any other $\rho_{O}$-cycles.

Let $B = \{b_{0}, \ldots, b_{k}\}$ be some class of equivalence of the relation $\hat\rho_{O}$ and the following correlations are satisfied
\[
b_{0} \rho b_{1},\, \ldots,\, b_{k-1} \rho b_{k},\, b_{k} \rho b_{0} \,.
\]
Support that $B \varsubsetneq A$. Let us fix $a \in A \setminus B$. By definition the following correlations hold true
\begin{align*}
& O(b_{i-1}, b_{i}, a) \,, \quad i \in \{1, \ldots, k\} \,; \\
& O(b_{k}, b_{0}, a) \,.
\end{align*}
Thus it follows from definition of cyclic order it follows that
\begin{align*}
& O(a, b_{i-1}, b_{i}) \,, \quad i \in \{1, \ldots, k\} \,; \\
& O(a, b_{k}, b_{0}) \,.
\end{align*}

From definition it also follows that if $O(a, b_{0}, b_{i-1})$ and $O(a, b_{i-1}, b_{i})$, then $O(a, b_{0}, b_{i})$. Therefore starting from $O(a, b_{0}, b_{1})$ in the finite number of steps we get $O(a, b_{0}, b_{k})$.

Thus $O(a, b_{k}, b_{0})$ and $O(a, b_{0}, b_{k})$ should be satisfied simultaneously but it contradicts to antisymmetry of cyclic order.

Therefore all elements of a set $A$ are equivalent under $\hat\rho_{O}$ and generate $\rho_{O}$-cycle.
\end{proof}


\begin{ozn}\label{ozn_induc_by_zruchn}
Let $\rho$ be a convenient relation on a finite set $A$. We define a ternary relation $O_{\rho}$ on $A$ with the help of the following rule. The ordered triple $(a_{1}, a_{2}, a_{3})$ of $A$ is said to be in the relation $O_{\rho}$ if $a_{1} \neq a_{2} \neq a_{3} \neq a_{1}$ and there are
\begin{equation}\label{eq_induc_by_zruchn}
a_{1} = a^{12}_{0}, a^{12}_{1}, \ldots, a^{12}_{m(1)} = a_{2} = a^{23}_{0}, \ldots, a^{23}_{m(2)} = a_{3} = a^{31}_{0}, \ldots, a^{31}_{m(3)} = a_{1} \,,
\end{equation}
which satisfy the following conditions:
\begin{itemize}
	\item $a^{sr}_{n-1} \rho a^{sr}_{n}$ for all $n \in \{1, \ldots, m(s)\}$ and $(s+1) \equiv r \pmod{3}$;
	\item $a^{sr}_{n} \notin \{a_{1}, a_{2}, a_{3}\}$ for all $n \in \{1, \ldots, m(s)-1\}$ and $(s+1) \equiv r \pmod{3}$.
\end{itemize}
\end{ozn}

\begin{prop}\label{prop_induced_cyclic}
The relation $O_{\rho}$ is a cyclic order on $A$.
\end{prop}

\begin{proof}
From definition it is obvious that the relation $O_{\rho}$ is cyclic.

Let us remark that from definition if $O_{\rho}(a_{1}, a_{2}, a_{3})$, then all elements of a set~\eqref{eq_induc_by_zruchn} (in particular elements $a_{1}$, $a_{2}$ and $a_{3}$) belong to the same class of equivalence of minimal equivalence relation $\hat\rho$ which contains $\rho$.

We should verify that all elements $a^{sr}_{n}$, $n \in \{1, \ldots, m(s)\}$, $(s+1) \equiv r \pmod{3}$ are different.

Suppose that it is not true and there are two different sets of indexes such that $a^{sr}_{n} = a^{t\tau}_{k}$, $n \in \{1, \ldots, m(s)\}$, $k \in \{1, \ldots, m(t)\}$, $(s+1) \equiv r \pmod{3}$, $(t+1) \equiv \tau \pmod{3}$.

Let us consider two sequences
\begin{align*}
(b_{1}, \ldots, b_{i}) & = (a^{sr}_{n}, a^{sr}_{n+1}, \ldots, a^{sr}_{m(s)}, \ldots, a^{t\tau}_{0}, a^{t\tau}_{1}, \ldots, a^{t\tau}_{k-1}, a^{t\tau}_{k}) \,, \\
(c_{1}, \ldots, c_{j}) & = (a^{t\tau}_{k}, a^{t\tau}_{k+1}, \ldots, a^{t\tau}_{m(t)}, \ldots, a^{sr}_{0}, a^{sr}_{1}, \ldots, a^{sr}_{n-1}, a^{sr}_{n}) \,.
\end{align*}
Those two sequences satisfy the following conditions:
\begin{itemize}
	\item $b_{l-1} \rho b_{l}$ for all $l \in \{1, \ldots, i\}$;
	\item $c_{l-1} \rho c_{l}$ for all $l \in \{1, \ldots, j\}$;
	\item $b_{i} = c_{1} = c_{j} = b_{1}$;
	\item there exists $\hat{a} \in \{a_{1}, a_{2}, a_{3}\}$ such that either $\hat{a} \in \{b_{1}, \ldots, b_{i}\} \setminus \{c_{1}, \ldots, c_{j}\}$ or $\hat{a} \in \{c_{1}, \ldots, c_{j}\} \setminus \{b_{1}, \ldots, b_{i}\}$  since  every element $a_{1}$, $a_{2}$, $a_{3}$ is contained exactly once in the sequence~\eqref{eq_induc_by_zruchn} by definition.
\end{itemize}

Let $\hat{a} \notin \{b_{1}, \ldots, b_{i}\}$. By definition the elements $b_{1}, \ldots, b_{i}$ generate a cycle therefore the set $\{b_{1}, \ldots, b_{i}\}$ is a class of equivalence of the relation $\hat\rho$, see Corollary~\ref{nas_struct_2}. But it contradicts to the condition that all elements of the set~\eqref{eq_induc_by_zruchn} belong to the same class of equivalence of the relation $\hat\rho$.

The case when $\hat{a} \notin \{c_{1}, \ldots, c_{j}\}$ can be considered similarly.

Therefore all elements of the set~\eqref{eq_induc_by_zruchn} are different.

Let $O_{\rho}(a_{1}, a_{2}, a_{3})$ and $O_{\rho}(a_{3}, a_{2}, a_{1})$ hold true simultaneously.
Then from definition there are  two sequences $a_{3} = a^{31}_{0}, a^{31}_{1}, \ldots, a^{31}_{m(3)} = a_{1}$ and $a_{1} = b^{31}_{0}, b^{31}_{1}, \ldots, b^{31}_{n(3)} = a_{3}$ such that
\begin{itemize}
	\item $a^{31}_{i-1} \rho a^{31}_{i}$ for all $i \in \{1, \ldots, m(3)\}$;
	\item $b^{31}_{j-1} \rho b^{31}_{j}$ for all $j \in \{1, \ldots, n(3)\}$;
	\item $a_{2} \notin \{a^{31}_{0}, \ldots, a^{31}_{m(3)}, b^{31}_{0}, \ldots, b^{31}_{n(3)}\}$.
\end{itemize}
It is obvious that there is $k \in \{1, \ldots, m(3)\}$ such that $a^{31}_{i} \notin \{b^{31}_{0}, \ldots, b^{31}_{n(3)}\}$ for $i < k$ but $a^{31}_{k} \in \{b^{31}_{0}, \ldots, b^{31}_{n(3)}\}$. Hence $a^{31}_{k} = b^{31}_{l}$ for some $l \in \{1, \ldots, n(3)\}$ and $a^{31}_{k} \rho b^{31}_{l+1}$. It is clear that all elements of the following sequence
\[
a_{3} = a^{31}_{0}, \ldots, a^{31}_{k}, b^{31}_{l+1}, \ldots, b^{31}_{n(3)}
\]
are different and generate $\rho$-cycle. Further by definition $a_{2}$ does not belong to that sequence. Therefore $a_{3} = a^{31}_{0}$ and $a_{2}$ belong to different classes of equivalence of relation $\hat\rho$, see Corollary~\ref{nas_struct_2}.

On the other hand elements $a_{1}$, $a_{2}$ and $a_{3}$ must belong to the unique class of equivalence $\hat\rho$, see above.

This contradiction proves the antisymmetry of the relation $O_{\rho}$.

Let $O_{\rho}(a_{1}, a_{2}, a_{3})$ and $O_{\rho}(a_{1}, a_{3}, a_{4})$ for some $a_{1}, \ldots, a_{4} \in A$.

We should remark that the elements $a_{1}, \ldots, a_{4}$ are pairwise different. Really, by definition $a_{1} \neq a_{3}$ and $\{a_{1}, a_{3}\} \cap \{a_{2}, a_{4}\} = \emptyset$. If $a_{2} = a_{4}$, then from a cyclicity of relation $O_{\rho}$ it follows that $O_{\rho}(a_{3}, a_{1}, a_{2})$ and $O_{\rho}(a_{4}, a_{1}, a_{3}) = O_{\rho}(a_{2}, a_{1}, a_{3})$. But it is impossible since a relation $O_{\rho}$ is antisymmetric.

Let us consider a sequence~\eqref{eq_induc_by_zruchn}. Its elements generate $\rho$-cycle. We will prove that $a_{4} \in \{a^{31}_{1}, \ldots, a^{31}_{m(3)-1}\}$.

Suppose that $a_{4} \in \{a^{12}_{1}, \ldots, a^{12}_{m(1)-1}\}$. Then $a_{4} = a^{12}_{k}$, $k \in \{1, \ldots, m(1)-1\}$. We consider the sequences
\begin{align*}
(b^{12}_{0}, \ldots, b^{12}_{t(1)}) & = (a_{1} = a^{12}_{0}, \ldots, a^{12}_{k} = a_{4}) \,; \\
(b^{23}_{0}, \ldots, b^{23}_{t(2)}) & = (a_{4} = a^{12}_{k}, \ldots, a^{12}_{m(1)} = a^{23}_{0}, \ldots, a^{23}_{m(2)} = a_{3}) \,;\\
(b^{31}_{0}, \ldots, b^{31}_{t(3)}) & = (a_{3} = a^{31}_{0}, \ldots, a^{31}_{m(3)} = a_{1}) \,.
\end{align*}

Join them into a sequence
\[
a_{1} = b^{12}_{0}, \ldots, b^{12}_{t(1)} = a_{4} = b^{23}_{0}, \ldots, b^{23}_{t(2)} = a_{3} = b^{31}_{0}, \ldots, b^{31}_{t(3)} \,.
\]
By the construction all elements of such sequence generate $\rho$-cycle therefore it satisfies the properties which are similar to the conditions of the sequence~\eqref{eq_induc_by_zruchn}. We get $O_{\rho}(a_{1}, a_{4}, a_{3})$. Then from a cyclicity of the relation $O_{\rho}$ it follows that  $O_{\rho}(a_{4}, a_{3}, a_{1})$. But by the condition we have $O_{\rho}(a_{1}, a_{3}, a_{4})$, moreover, we proved that the relation $O_{\rho}$ is antisymmetric. Thus the relation $O_{\rho}(a_{4}, a_{3}, a_{1})$ does not hold true and $a_{4} \notin \{a^{12}_{1}, \ldots, a^{12}_{m(1)-1}\}$.

The fact that $a_{4} \notin \{a^{23}_{1}, \ldots, a^{23}_{m(2)-1}\}$ can be proved similarly.

Therefore $a_{4} \in \{a^{31}_{1}, \ldots, a^{31}_{m(3)-1}\}$ and $a_{4} = a^{31}_{s}$ for some $s \in \{1, \ldots, m(3)-1\}$.

Let us consider the sequences
\begin{align*}
(c^{12}_{0}, \ldots, c^{12}_{\tau(1)}) & = (a_{1} = a^{12}_{0}, \ldots, a^{12}_{m(1)} = a_{2}) \,; \\
(c^{23}_{0}, \ldots, c^{23}_{\tau(2)}) & = (a_{2} = a^{23}_{0}, \ldots, a^{23}_{m(2)} = a^{31}_{0}, \ldots, a^{31}_{s} = a_{4} ) \,;\\
(c^{31}_{0}, \ldots, c^{31}_{\tau(3)}) & = (a_{4} = a^{31}_{s}, \ldots, a^{31}_{m(3)} = a_{1}) \,.
\end{align*}
Let us join them into a sequence
\[
a_{1} = c^{12}_{0}, \ldots, c^{12}_{\tau(1)} = a_{2} = c^{23}_{0}, \ldots, c^{23}_{\tau(2)} = a_{4} = c^{31}_{0}, \ldots, c^{31}_{\tau(3)} \,.
\]
By construction this sequence satisfies the conditions of definition~\ref{ozn_induc_by_zruchn}. Therefore the correlation $O_{\rho}(a_{1}, a_{2}, a_{4})$ holds true and the relation $O_{\rho}$ is transitive.

Finally, we can conclude that the relation $O_{\rho}$ satisfies all conditions of definition of cyclic order.
\end{proof}

\begin{ozn}\label{ozn_cyclic_morphisms}
Let $C$ and $D$ be cyclic orders on sets $A$ and $B$, respectively. Let $\varphi : A \rightarrow B$ be a bijective map.

A map $\varphi$ is called \emph{a monomorphism} of cyclic order $C$ into a cyclic order $D$ if $C(a_{1}, a_{2}, a_{3}) \Rightarrow D(\varphi(a_{1}), \varphi(a_{2}), \varphi(a_{3}))$; it is  called \emph{an epimorphism } $C$ onto $D$ if $D(b_{1}, b_{2}, b_{3}) \Rightarrow C(\varphi^{-1}(b_{1}), \varphi^{-1}(b_{2}), \varphi^{-1}(b_{3}))$; $\varphi$ is \emph{an isomorphism} $C$ onto $D$ if $C(a_{1}, a_{2}, a_{3}) \Leftrightarrow \\D(\varphi(a_{1}), \varphi(a_{2}), \varphi(a_{3}))$.
\end{ozn}

\begin{remk}
It is clear that

1) if $\varphi$ is a monomorphism of cyclic order $C$ onto $D$, then $\varphi^{-1}$ is an epimorphism of $D$ onto $C$;

2) an isomorphism of the relations of cyclic order is a map which is a monomorphism and an epimorphism simultaneously;

3) a relation of isomorphism is a relation of equivalence.
\end{remk}

\begin{lem}\label{lem_full_cycle_iso}
Let $C$ and $D$ be complete cyclic orders on the sets $A$ and $B$, respectively, $\varphi : A \rightarrow B$ is a bijective map.

If $\varphi$ is either monomorphism or an epimorphism, then $\varphi$ is an isomorphism.
\end{lem}

\begin{proof}
Let $\varphi$ be an epimorphism (in the case when $\varphi$ is a monomorphism we consider a map $\varphi^{-1}$). Let us check that $\varphi$ is also a monomorphism.

Let $C(a_{1}, a_{2}, a_{3})$ for some $a_{1}$, $a_{2}$, $a_{3} \in A$. We define $b_{i} = \varphi(a_{i}) \in B$, $i = 1, 2, 3$. From definition it follows that $a_{1} \neq a_{2} \neq a_{3} \neq a_{1}$. Then $b_{1} \neq b_{2} \neq b_{3} \neq b_{1}$.

The cyclic order $D$ is full therefore there is a permutation $\sigma \in S(3)$ such that $D(b_{\sigma(1)}, b_{\sigma(2)}, b_{\sigma(3)})$. From an epimorphism of $\varphi$ we can conclude that $C(a_{\sigma(1)}, a_{\sigma(2)}, a_{\sigma(3)})$. Thus $\sigma$ is even permutation. Now from antisymmetry an cyclicity of $D$ it follows that $D(b_{1}, b_{2}, b_{3})$, see\cite{Nov}. Therefore we get $D(\varphi(a_{1}),\varphi(a_{2}), \varphi(a_{3}))$ and $\varphi$ is a monomorphism.
\end{proof}

\begin{remk}
Lemma~\ref{lem_full_cycle_iso} holds true for arbitrary sets $A$ and $B$, i.e. they can be infinite.
\end{remk}

\begin{lem}\label{lem_induced_induced_equiv}Let $O$ is a relation of complete cyclic order on the finite set $A$. Then
\[
O = O_{\rho_{O}} \,,
\]
 where $\rho_{O}$ is a convenient binary relation generated by $O$ and $O_{\rho_{O}}$ is a relation of cyclic order generated by the convenient relation $\rho_{O}$.
\end{lem}

\begin{proof}
We should prove that the relation $O_{\rho_{O}}$ is full.

Let $b_{1}$, $b_{2}$, $b_{3}$ be some pairwise different elements of $A$. From Proposition~\ref{prop_cycle} and Corollary~\ref{nas_struct_2} the minimal relation of equivalence $\hat\rho_{O}$ which contains $\rho_{O}$ has the unique class of equivalence $B = A$. Thus we can index all elements of $A$ in such way that~\eqref{eq_max_chain} holds true. From Corollary~\ref{nas_struct_1} we also get $a_{n} \rho a_{0}$.

It is obvious that $\{b_{1}, b_{2}, b_{3}\} = \{a_{k_{1}}, a_{k_{2}}, a_{k_{3}}\}$ for some $0 \leq k_{1} < k_{2} < k_{3} \leq n$ further there is a inversion $\sigma \in S(3)$ such that $a_{k_{i}} = b_{\sigma(i)}$, $i=1, 2, 3$.

Let us consider the sequences
\begin{align*}
(c^{12}_{0}, \ldots, c^{12}_{m(1)}) & = (a_{k_{1}}, a_{k_{1}+1}, \ldots, a_{k_{2}}) \,; \\
(c^{23}_{0}, \ldots, c^{23}_{m(2)}) & = (a_{k_{2}}, \ldots, a_{k_{3}}) \,; \\
(c^{31}_{0}, \ldots, c^{31}_{m(3)}) & = (a_{k_{3}}, \ldots, a_{n}, a_{0}, \ldots, a_{k_{1}}) \,.
\end{align*}
We can join them into one
\[
a_{k_{1}} = c^{12}_{0}, \ldots, c^{12}_{m(1)} = a_{k_{2}} = c^{23}_{0}, \ldots, c^{23}_{m(2)} = a_{k_{3}} = c^{31}_{0}, \ldots, c^{31}_{m(3)} = a_{k_{1}} \,.
\]
By construction this sequence satisfies the conditions of Definition~\ref{ozn_induc_by_zruchn} thus we get $O_{\rho_{O}}(a_{k_{1}}, a_{k_{2}}, a_{k_{3}})$. It means that $O_{\rho_{O}}(b_{\sigma(1)},b_{\sigma(2)}, b_{\sigma(3)})$ and $O_{\rho_{O}}$ is full.

Suppose that $O_{\rho_{O}}(a_{1}, a_{2}, a_{3})$ holds true for some $a_{1}$, $a_{2}$, $a_{3} \in A$. From Definition~\ref{ozn_induc_by_zruchn} it follows that there is a sequence
\[
a_{1} = a^{12}_{0}, \ldots, a^{12}_{m(1)} = a_{2} \,,
\]
  such that $a^{12}_{i-1} \rho_{O} a^{12}_{i}$ for all $i \in \{1, \ldots, m(1)\}$. Therefore from definition of the relation $\rho_{O}$  correlations $O(a^{12}_{i-1}, a^{12}_{i}, a)$ follow for all $a \in A \setminus \{a^{12}_{i-1}, a^{12}_{i}\}$, $i \in \{1, \ldots, m(1)\}$. In particular,  $O(a^{12}_{i-1}, a^{12}_{i}, a_{3})$, $i \in \{1, \ldots, m(1)\}$. From cyclicity of $O$ it follows that the correlations $O(a_{3}, a^{12}_{i-1}, a^{12}_{i})$, $i \in \{1, \ldots, m(1)\}$ hold true.

Starting from the correlation $O(a_{3}, a^{12}_{0}, a^{12}_{1}) = O(a_{3}, a_{1}, a^{12}_{1})$, using the previous correlations and transitivity of $O$ we inductively get that $O(a_{3}, a_{1}, a^{12}_{i})$, $i \in \{1, \ldots, m(1)\}$. In particular, $O(a_{3}, a_{1}, a^{12}_{m(1)}) = O(a_{3}, a_{1}, a_{2})$. From a cyclicity of $O$ it follows that $O(a_{1}, a_{2}, a_{3})$.

Therefore an identical map $Id_{A} : A \rightarrow A$ induces an epimorphism of a complete cyclic order $O$ onto a complete cyclic order $O_{\rho_{O}}$. From Lemma~\ref{lem_full_cycle_iso} it follows that the map $Id_{A}$ is an isomorphism of the cyclic orders $O$ and $O_{\rho_{O}}$ therefore $O = O_{\rho_{O}}$.
\end{proof}


\begin{lem}\label{lem_two_pairs}
Let $\rho$ be a convenient relation such that all elements of a set $A$, $\sharp A \geq 3$ generate a cycle.

Suppose that a graph of relation $\mu$ on $A$ is obtained from a graph of $\rho$ by throwing out two pairs $(b_{1}, b_{1}')$ and $(b_{2}, b_{2}')$ (the cases when either $b_{1}' = b_{2}$ or $b_{2}' = b_{1}$ are included). Let $\hat\mu$ be a minimal relation of equivalence which contains $\mu$.

Then the relation $\mu$ is convenient, $\hat\mu$ has exactly two classes of equivalence $B_{1}$ and $B_{2}$ such that the elements of each of them generate $\mu$-chain and the elements $b_{1}$, $b_{2} \in A$ belong to the different classes of equivalence of $\hat\mu$.
\end{lem}

\begin{proof}
The fact that $\mu$ is a convenient relation is trivial corollary from definition.

The relation $\mu$ does not contain cycles. In fact, if the elements of some set $B \subseteq A$ generate $\mu$-cycle, then elements of $B$ generate $\rho$-cycle. From Corollary~\ref{nas_struct_2} and the condition of lemma we get $B = A$. Then from definition of a cycle it follows that there is $a \in A$ such that $b_{1} \mu a$, hence $b_{1} \rho a$. But $b_{1} \rho b_{1}'$ and $b_{1}' \neq a$ (by condition of lemma $b_{1}$ is not in the relation $\mu$ with $b_{1}'$). It contradicts to the Condition 2) of definition~\ref{ozn_zruchne}.

Thus every class of equivalence of the relation $\hat\mu$ is a chain, see Corollary~\ref{nas_struct_1}, and it contains exactly one  element which is in the relation $\mu$ with no element of $A$.

By condition of lemma the elements of $A$ generate $\rho$-cycle. From Definition~\ref{ozn_zruchne} it follows that there is the unique $a' \in A$ such that $a \rho a'$ for every $a \in A$. Then $a \mu a'$, if $a \notin \{b_{1}, b_{2}\}$ but $b_{1}$ and $b_{2}$ are the unique elements of the set $A$ which are not in the relation $\mu$ with any element of $A$.

now the statement of lemma elementary follows from what we said before.
\end{proof}


\subsection{A local connectivity of two dimensional disk in boundary points}

\begin{ozn}~\cite{Newman, Z-F-C}
Let $E$ be a subset of a topological space $S$ and $x$ is some point of $S$ ($x$ does not necessarily belong to $E$). A set $E$ is called \emph{a locally connected in a point $x$} if for every neighborhood $U$ of $x$ there is a neighborhood $U' \subseteq U$ of $x$ such that any two points which belong to $U' \cap E$ can be joined by a connected set which belongs to $U \cap E$.
\end{ozn}

\begin{lem}\label{lem_loc_connect}
Let $D^{2}$ be a closed two dimensional disk, $x \in \partial D^{2}$ and $W$ an open neighborhood of point $x$ in a space $D^{2}$.

If for some connected components $W_{1}$ and $W_{2}$ of a set $W \cap (D^{2} \setminus \partial D^{2})$ the following correlation holds true $x \in \Cl{W}_{1} \cap \Cl{W}_{2}$, then $W_{1} = W_{2}$.
\end{lem}

\begin{proof}
Obviously, we can assume that $D^{2}$ is a standard two dimensional disk on a plane. Let $U$ be a neighborhood of point $x$ in $\rr^{2}$ such that  $D^{2} \cap U = W$. It is known, see~\cite{Newman, Z-F-C}, that every Jordan domain on the plane is locally connected in all points of its boundary. Therefore there exists a neighborhood $U'$ of $x$ such that arbitrary two points which belong to $U' \cap (D^{2} \setminus \partial D^{2})$ can be connected by a connected set that is contained in $U \cap (D^{2} \setminus \partial D^{2})$. Therefore all points of the set $U' \cap (D^{2} \setminus \partial D^{2})$ should belong to the unique connected component of a set $W \cap (D^{2} \setminus \partial D^{2})$.
\end{proof} 

\section{Criterion of a $\cD$-planarity of a tree.}

Let $T$ be a tree, $V$ a set of its vertices, $V_{ter}$ a set of its terminal vertices and $V^{\ast} \subseteq V$ a subset of $T$ such that $V_{ter} \subseteq V^{\ast}$. We assume that if $\sharp V^{\ast} \geq 3$ then there is some cyclic order $C$ defined on $V^{\ast}$.

Let
\[
D^{2} = \{ (x,y) \in \rr^{2} \,|\, x^{2} + y^{2} \leq 1 \}
\]
be a closed oriented 2--dimensional disk.

\begin{ozn}\label{ozn_01}
 A tree $T$ is called \emph{$\cD$-planar} if there exists an embedding $\varphi : T \rightarrow \rr^{2}$ which satisfies~\eqref{eq_02} and if $\sharp V^{\ast} \geq 3$ then a cyclic order $\varphi(C)$ on $\varphi(V^{\ast})$ coincides with a cyclic order which is generated by the orientation of $\partial D^{2} \cong S^{1}$.
\end{ozn}

\begin{remk}\label{rem_01}
 A map $\varphi|_{V^{\ast}} : V^{\ast} \rightarrow \varphi(V^{\ast})$ is bijective whence a ternary relation $\varphi(C)$ on $\varphi(V^{\ast})$ defined by following correlation
\[
C(v_{1}, v_{2}, v_{3}) \Rightarrow \varphi(C)(\varphi(v_{1}), \varphi(v_{2}), \varphi(v_{3})) \,, \quad
v_{1}, v_{2}, v_{3} \in V^{\ast} \,,
\]
is a relation of cyclic order.
\end{remk}

\begin{remk}\label{rem_02}
 We can define a cyclic order in a natural way on an oriented circle $S^{1}$: an ordered triple of points $x_{1}$, $x_{2}$, $x_{3} \in S^{1}$ is cyclically ordered if these points are passed in that order in the process of moving along a circle in a positive direction.
\end{remk}

\begin{figure}[htbp]
\centerline{\includegraphics{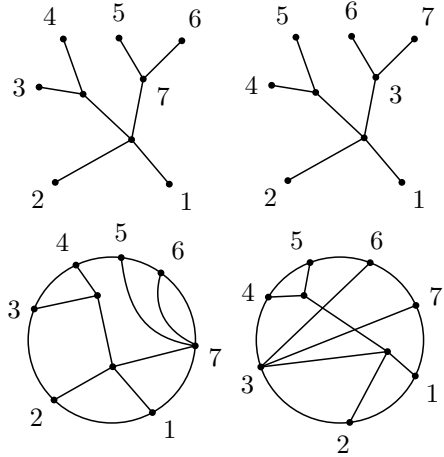}}
\begin{center}
\caption{On the left a tree is $\cD$-planar.}
 \end{center}
\end{figure}

\begin{theorem}\label{crit_d_planarity}
If $V^{\ast}$ contains just two vertices, a tree $T$ is $\cD$-planar.

If $\sharp V^{\ast} \geq 3$ then a $\cD$-planarity of $T$ is equivalent to satisfying the following condition:
\begin{itemize}
    \item for any edge $e$ there are exactly two paths such that they pass through an edge $e$ and connect two adjacent vertices of $V^{\ast}$.
\end{itemize}
\end{theorem}

\begin{proof}
 If $\sharp V^{\ast} = 2$, then $T$ is homeomorphic to a segment and a set of its terminal vertices coincides with $V^{\ast} = V_{ter}$, see Lemma~\ref{lemma_01}. It is obvious that there exists an embedding $\varphi : T \rightarrow \rr^{2}$ satisfying Definition~\ref{ozn_01} and a tree $T$ is $\cD$-planar.

 Let $\sharp V^{\ast} \geq 3$ and $T$ is $\cD$-planar. It means that there is an embedding $\varphi : T \rightarrow \rr^{2}$ which satisfies Definition ~\ref{ozn_01}.

Let $e \in E(T)$ be an edge of $T$ connecting vertices $w_{1}$, $w_{2} \in V$. We fix a point $x \in \varphi(e) \setminus \{\varphi(w_{1}), \varphi(w_{2})\}$.

A topological space $T$ is one--dimensional compact hence its homeomorphic image $\varphi(T)$ is one--dimensional~\cite{Kuratowsky_V1}. Then $x \in \Cl{(\rr^{2} \setminus \varphi(T))}$. It follows from~\eqref{eq_02} that $x \in \Int{D^{2}}$, therefore
\[
x \in \Cl{(\rr^{2} \setminus (\varphi(T) \cup \partial D^{2}))} \,.
\]
  By Lemma~\ref{lemma_01} there is a connected component $U_{j}$ of a set $\rr^{2} \setminus (\varphi(T) \cup \partial D^{2})$ such that a point $x$ belongs to a boundary of it.

Corollary~\ref{nas_01} states that $\partial U_{j} \cap \varphi(T) = \varphi(P(v_{j}, v_{j}'))$, where $\varphi(v_{j})$, $\varphi(v_{j}')$ are adjacent with respect to a cyclic order of $\varphi(V^{\ast})$  induced from $\partial D^{2}$, see Remark~\ref{rem_02}. According to Definition~\ref{ozn_01}, it is the same as vertices $v_{j}$ and $v_{j}'$ are adjacent under a cyclic order $C$ on $V^{\ast}$.

  So $e \in P(v_{j}, v_{j}')$ and vertices $v_{j}$, $v_{j}'$ are adjacent. It means that for any edge of a $\cD$-planar tree $T$ there is at least one path that satisfies the condition of theorem.

 There exists an open neighborhood $W = e \setminus \{w_{1}, w_{2}\}$ of a point $\varphi^{-1}(x)$ in $T$  that is homeomorphic to an interval. Using the compactness of $T \setminus W$ and theorem of Shenflies~\cite{Newman, Z-F-C} we can find a neighborhood $U$ of $x$ in $\rr^{2} \setminus \partial D^{2}$ and a homeomorphism $h : U \rightarrow \Int{D^{2}}$ such that $h(x) = (0,0)$, $h \circ \varphi(T) = h \circ \varphi(W) = (-1, 1) \times \{0\}$. Let us designate
\begin{gather*}
U^{+} = h^{-1}(\{(x,y) \in \Int{D^{2}} \,|\, y > 0\}) \,, \\
U^{-} = h^{-1}(\{(x,y) \in \Int{D^{2}} \,|\, y < 0\}) \,.
\end{gather*}

 It is clear that $U \subseteq \varphi(T) \cup U^{+} \cup U^{-}$. If for some component $U_{k}$ of $\rr^{2} \setminus (\varphi(T) \cup \partial D^{2})$ the intersections $U^{+} \cap U_{k}$ and $U^{-} \cap U_{k}$ are empty, then $x \notin \Cl{U}_{k}$ and $e \notin P(v_{k}, v_{k}')$ in terms of Lemma~\ref{lemma_01}.

By the construction, the sets  $U^{+}$ and $U^{-}$ are connected and they belong to $\rr^{2} \setminus (\varphi(T) \cup \partial D^{2})$. Thus there are two components $U_{i}$ and $U_{j}$ such that $U^{+} \in U_{i}$, $U^{-} \in U_{j}$, $x \in \Cl{U}_{i} \cap \Cl{U}_{j}$ and $e \in P(v_{i}, v_{i}') \cap P(v_{j}, v_{j}')$.

 By Corollaries~\ref{nas_01} and~\ref{nas_02} for any edge of $T$ there are no more then two paths such that they connect adjacent vertices of $V^{\ast}$.

In order to verify that there are exactly two such paths it is sufficient to prove that $U_{i} \neq U_{j}$.

Suppose that for some component $U_{i}$  of $\rr^{2} \setminus (\varphi(T) \cup \partial D^{2})$ we get $U \setminus \varphi(T) = U^{+} \cup U^{-} \subseteq U_{i}$. An open connected subset $U_{i}$ of  $\rr^{2}$ is path-connected~\cite{Kuratowsky_V2}.

 Denote $a_{0}^{+} = (0, 1/2)$, $a_{0}^{-} = (0, -1/2) \in \Int{D^{2}}$, $\gamma_{0} = \{0\} \times [-1/2, 1/2] \subseteq \Int{D^{2}}$, $a^{+} = h^{-1}(a_{0}^{+})$, $a^{-} = h^{-1}(a_{0}^{-}) \in U$, $\gamma = h^{-1}(\gamma_{0})$.

 It is obvious that the points $a_{0}^{+}$ and $a_{0}^{-}$ are attainable from domains $h(U^{+}) \setminus \gamma_{0}$ and $h(U^{-}) \setminus \gamma_{0}$ by a simple continuous curve. Therefore the points $a^{+}$ and $a^{-}$ are attainable from the domain $U_{i} \setminus \gamma$ and there is a cut $\hat{\gamma}$ of $U_{i} \setminus \gamma$ between $a^{+}$ and $a^{-}$~\cite{Newman, Z-F-C}.

 Then $\mu = \gamma \cup \hat{\gamma}$ is a simple close curve such that $\mu \cap \varphi(T) = \{x\}$, $\mu \setminus \{x\} \subseteq U_{i}$ and $h(\mu) \supseteq \gamma_{0}$.

By Jordan's theorem $\mu$ bounds an open disk $G$~\cite{Newman, Z-F-C}.

The point $x$ does not belong to the compact $\hat{\gamma}$ hence there exists its open neighborhood $\hat{U} \subseteq U$ such that $\hat{U} \cap \hat{\gamma} = \emptyset$. Since $h$ maps a neighborhood $\hat{U}$ of a point $x$ into an open neighborhood of origin then there exists an $\varepsilon \in (0, 1/2)$ such that a set
\[
Q_{0} = \{(x,y) \in D^{2} \,|\, x^{2}+y^{2} < \varepsilon^{2}\}
\]
does not intersect the set $h(\hat{\gamma})$. It follows that
\begin{gather*}
Q_{0} \cap h(\varphi(T) \cup \partial D^{2}) = Q_{0} \cap h \circ \varphi(e) = (-\varepsilon, \varepsilon) \times \{0\} \,, \\
Q_{0} \cap h(\mu) = Q_{0} \cap \gamma_{0} = \{0\} \times (-\varepsilon, \varepsilon) \,.
\end{gather*}
Denote $Q = h^{-1}(Q_{0})$. Evidently, a set $Q$ is an open neighborhood of $x$.

Open sets
\[
h^{-1} \left( \{(x,y) \in Q_{0} \,|\, x < 0\} \right)
\mbox{ and }
h^{-1} \left( \{(x,y) \in Q_{0} \,|\, x > 0\} \right)
\]
are connected and do not intersect the set $\mu$. Therefore one of them must be contained in a disk $G$, another should belong to an unbounded domain $\rr^{2} \setminus \Cl{G}$.

Sets $h^{-1}((-\varepsilon, 0) \times \{0\})$ and $h^{-1}((0, \varepsilon) \times \{0\})$ belong to the intersection of these domains with the image $\varphi(e)$ of $e$. Hence $\varphi(e) \cap G \neq \emptyset$ and $\varphi(e) \cap \rr^{2} \setminus \Cl{G} \neq \emptyset$ hold true.

A segment $\varphi(e)$ is divided by $x$ on two connected arcs that have no common points with $\mu = \partial G$. Thus one of them should belong to $G$ and the other is contained in $\rr^{2} \setminus \Cl{G}$.

Finally, the following statement is true: either  $\varphi(w_{1})$ or $\varphi(w_{2})$ belongs to $G$ and the other point is contained in $\rr^{2} \setminus \Cl{G}$.

Let $\varphi(w_{1}) \in G$, $\varphi(w_{2}) \in \rr^{2} \setminus \Cl{G}$.

  By the construction, curves $\partial D^{2}$ and $\mu$ have no common points since either $G \subseteq \Int{D^{2}}$ or $\Int{D^{2}} \subseteq G$. But $\emptyset \neq (\gamma \cap Q) \subseteq (\mu \cap U_{i}) \subseteq (\mu \cap \Int{D^{2}})$. Therefore $\Cl{G} \subseteq \Int{D^{2}}$.

Let us denote by  $\hat{T}$ a graph with a set of vertices $V(\hat{T}) = V(T) = V$ and a set of edges $E(\hat{T}) = E(T) \setminus \{e\} = E \setminus \{e\}$.

 It is easy to show that the graph $\hat{T}$ has two connected components $T_{1} \ni w_{1}$ and $T_{2} \ni w_{2}$. The images of them do not intersect with the curve $\mu$, therefore a set $\varphi(T_{1})$ together with the point $\varphi(w_{1})$ belongs to $G \subseteq \Int{D^{2}}$ and $\varphi(T_{2}) \subseteq \rr^{2} \setminus \Cl{G}$.

By relation $\varphi(w_{1}) \in G \subseteq \Int{D^{2}}$ and Condition~\eqref{eq_02}, the vertex $w_{1}$ has degree at least 2. Therefore it is adjacent to at least one edge of $T$ except $e$ that is an edge of $T_{1}$. It means that a tree $T_{1}$ is non degenerated.

Since degrees of all other vertices of $T_{1}$ in $T$ are the same as degrees in $T_{1}$ then $V_{ter}(T_{1}) \subseteq V_{ter}(T) \cup \{w_{1}\}$. As we know $\sharp V_{ter}(T_{1}) \geq 2$ whence there is $w \in V_{ter}(T_{1}) \cap V_{ter}(T)$.

By the construction $\varphi(w) \in G \subseteq Int{D^{2}}$.

On the other hand it follows from ~\eqref{eq_01} and~\eqref{eq_02} that $\varphi(w) \in \varphi(V^{\ast}) \subseteq \partial D^{2}$.

We have  the contradiction with the assumption that $U \setminus \varphi(T) \subseteq U_{i}$ for some $i$.

 So, there are exactly two components $U_{i} \neq U_{j}$ of a compliment $\rr^{2} \setminus (\varphi(T) \cup \partial D^{2})$ such that the point $x \in \varphi(e) \setminus \{\varphi(w_{1}), \varphi(w_{2})\}$ which is contained in the image of an edge $e$ of $T$ is a boundary point of. Consequently, by Corollaries~\ref{nas_01} and~\ref{nas_02} there are exactly two paths such that they pass through an arbitrary edge of $T$ and connect the adjacent vertices of $V^{\ast}$.

\medskip

  Let $\sharp V^{\ast} \geq 3$ and for any $e \in E(T)$ of $T$ there are exactly two paths such that they pass through this edge and connect adjacent vertices of $V^{\ast}$.

We should prove that the tree $T$ is $\cD$-planar.

At first we consider a relation $C$ that is a full cyclic order on a set $V^{\ast}$. It generates a convenient relation $\rho_{C}$ on  $V^{\ast}$, see Definition~\ref{ozn_induc_by_cycle}. Let us examine a set of the directed paths
\[
\cP = \{ P(v, v') \,|\, v' \rho_{C} v \}
\]
in $T$.

 By Definitions~\ref{ozn_cyclic_susidni} and~\ref{ozn_induc_by_cycle} two vertices $v$, $v' \in V^{\ast}$ are adjacent with respect to a cyclic order $C$  iff either $v \rho_{C} v'$ or $v' \rho_{C} v$ is true. These correlations can not hold true simultaneously,  since a pair of vertices $v$, $v'$ would generate a $\rho_{C}$-cycle, see Definition~\ref{ozn_cycle}, and this contradicts to Proposition~\ref{prop_cycle} and Corollary~\ref{nas_struct_2} since $\sharp V^{\ast} \geq 3$.

It follows from the discussion above that for every edge $e$ of $T$ there are exactly two paths of the set $\cP$ passing through $e$.

Let us consider a binary relation $\rho$ on the set $V_{\ast}$ which is defined by a correlation
\begin{equation}\label{eq_rho_by_paths}
v \rho v' \Leftrightarrow P(v, v') \in \cP \,.
\end{equation}

Evidently, relation $\rho$ is \emph{dual} to the relation $\rho_{C}$ ( $v_{1} \rho v_{2} \Leftrightarrow v_{2} \rho_{C} v_{1}$). Therefore by  Definition~\ref{ozn_zruchne},  $\rho$ is the convenient relation on $V^{\ast}$. So, a minimal relation of equivalence $\hat\rho$ on $V^{\ast}$ containing $\rho$ coincides with a  minimal relation of equivalence $\hat{\rho}_{C}$ on $V^{\ast}$ containing $\rho_{C}$. Thus the elements of the set $V^{\ast}$ generate a $\rho$-cycle, see Proposition~\ref{prop_cycle} and Corollary~\ref{nas_struct_1}.

 Let $e$ be an edge of the tree $T$. We should prove that those two directed paths of the set $\cP$ that contain $e$ pass through $e$ in opposite directions.

Let us consider a binary relation $\mu_{e}$ on $V^{\ast}$ that is defined as follows
\[
v \mu_{e} v' \Leftrightarrow P(v, v') \in \cP \text{ і } e \notin P(v, v') \,.
\]
It is easy to see that a diagram of the relation $\mu_{e}$ can be obtained from a diagram of $\rho$ by removing two pairs of vertices of $V^{\ast}$ corresponding to paths of $\cP$ which pass through $e$. Let $(v_{1}, v_{1}')$ and $(v_{2}, v_{2}')$ be such pairs. Therefore the relation $\mu_{e}$ satisfies the conditions of Lemma~\ref{lem_two_pairs}.

By this Lemma a minimal relation of equivalence $\hat{\mu}_{e}$ containing $\mu_{e}$ has two classes of equivalence $B_{1}$, $B_{2}$ and $v_{1} \in B_{1}$, $v_{2} \in B_{2}$.

 Let $w$, $w' \in V$ be the ends of $e$. Let us consider a subgraph $T'$ of the tree $T$ such that $V(T') = V(T)$ and $E(T') = E(T) \setminus \{e\}$. It is clear that the vertices $w$ and $w'$ belong to different connected components of a graph $T'$ (if there exists a path $P$ in $T'$ such that it connects them then these vertices can be connected by two different paths $P$ and $P' = \{e\}$ in the tree $T$ ). We denote these components by $T_{w}$ and $T_{w'}$.

Suppose that for vertices $v$, $v' \in V$ there is an directed path $P(v, v')$ passing through $e$. Let it first passes through the vertex $w$ and then though $w'$. Then paths $P(v, w)$ and $P(w', v')$ belong to $T'$, so $v \in T_{w}$, $v' \in T_{w'}$. In case when the path $P(v, v')$ first passes through $w'$ and then through $w$ we have $v' \in T_{w}$ and $v \in T_{w'}$.

 It is easy to see that every class of equivalence of the relation $\hat{\mu}_{e}$ belongs to the unique connected component of the set $T'$. By the construction different classes of equivalence have to belong to the different components of $T'$.

So, we conclude that either $B_{1} \subseteq T_{w}$ and $B_{2} \subseteq T_{w'}$ or $B_{1} \subseteq T_{w'}$ and $B_{2} \subseteq T_{w}$. Suppose that first pair of inequalities holds true.

If the directed paths $P(v_{1}, v_{1}')$ and $P(v_{2}, v_{2}')$ pass through $e$ in the same direction, then $P(v_{1}, w) \cup P(v_{2}, w) \subseteq T_{w}$ and $v_{2} \in T_{w}$. By the construction $T_{w} \cap V^{\ast} = B_{1}$ thus $v_{2} \in B_{1}$. But it is a contradiction to Lemma~\ref{lem_two_pairs}. So, the paths $P(v_{1}, v_{1}')$ and $P(v_{2}, v_{2}')$ pass through $e$ in the opposite directions.

The case $B_{1} \subseteq T_{w'}$, $B_{2} \subseteq T_{w}$ is considered similarly.

Let us construct an embedding of $T$ into oriented disk $D^{2}$.

 Let $D^{2}$ be an oriented disk (closed disk with a fixed orientation on the boundary), $I = [0, 1]$ an directed segment and $\psi : I \rightarrow D^{2}$ an embedding such that $\psi(I) \subseteq \partial D^{2}$. The direction of a segment is said to be \emph{coordinated} with the orientation of disk if a direction of passing along the simple continuous curve $\psi(I)$ from the origin $\psi(0)$ to the end $\psi(1)$ coincides with given orientation of the boundary $\partial D^{2}$.

Every directed path in $T$ is topologically a closed segment thus for directed path $P(v, v')$ with the origin $v$ and the end $v'$ there exists an embedding $\Phi_{P(v, v')} : P(v, v') \rightarrow D^{2}$ such that $\Phi_{P(v, v')}(P(v, v')) \subseteq \partial D^{2}$ and a direction of $P(v, v')$ is coordinated with the orientation of $D^{2}$.

We fix a disjoint union of closed oriented disks $\bigsqcup_{P \in \cP} D_{P}$ and a set of the embeddings
\begin{align}
\Phi_{P(v, v')} : P(v, v') & \rightarrow D_{P(v, v')} \,, \\
\Phi_{P(v, v')}(P(v, v')) & \subseteq \partial D_{P(v, v')} \,, \quad P(v, v') \in \cP \,, \nonumber
\end{align}
such that the directions of paths $P(v, v') \in \cP$ are coordinated with the orientations of corresponding disks.

Let  us consider a space
\[
\tilde{D} = T \sqcup \bigsqcup_{P \in \cP} D_{P} \,.
\]
All maps $\Phi_{P}$, $P \in \cP$ are injective therefore  a family of sets
\begin{equation*}
F_{x} = \left\{
\begin{array}{ll}
\{x\} \bigcup\limits_{P \in \cP \,:\, x \in P} \Phi_{P}(x) \,, & x \in T \,, \\
\{x\} \,, & x \in \bigcup\limits_{P \in \cP} D_{P} \setminus \Phi_{P}(P) \,.
\end{array}
\right.
\end{equation*}
generates a partition $\ff$ of the space $\tilde{D}$.

We consider a factor-space $D$ of $\tilde{D}$ over partition $\ff$ and a projection map
\[
\pr : \tilde{D} \rightarrow D \,.
\]

Let us prove that $D$ is homeomorphic to a disk, the orientations of disks $D_{P}$, $P \in \cP$ give some orientation on $D$ and a map
\[
\varphi = \restrict{\pr}{T} : T \rightarrow D
\]
conforms to the conditions of Definition~\ref{ozn_01}.

At first we investigate some properties of the space $D$ and the projection $\pr$.

\subsection{}\label{vlast_zamkn_pr} The mapping $\pr$ is closed.

Recall that a set is called \emph{saturated} over partition $\ff$ if it consists of entire elements of that partition.

Topology of space $D$ is a factor-topology (a set $A$ is closed in $D$ iff its full preimage $\pr^{-1}(A)$ is closed in $\tilde{D}$). For proof of closure of a projection map $\pr$ it is sufficient to check that for any closed subset $K$ of the space $\tilde{D}$ minimal saturated set $\tilde{K} = \pr^{-1}(\pr(K))$ containing $K$ is also closed.

From the definition of partition $\ff$ it follows that
\begin{eqnarray}\label{eq_preimage}
K & = & (K \cap T) \sqcup \bigsqcup_{P \in \cP} (K \cap D_{P}) \,, \nonumber \\
\tilde{K} & = & (K \cap T) \sqcup \bigsqcup_{P \in \cP} \left( (K \cap D_{P}) \cup \Phi_{P}(K \cap P) \right) \,.
\end{eqnarray}
 Sets $T$, $P$, $D_{P}$, $P \in \cP$ are compacts and all maps $\Phi_{P}$ are homeomorphisms onto their images. Thus all sets $K \cap T$, $K \cap D_{P}$, $\Phi_{P}(K \cap P)$, $P \in \cP$, are compacts. The graph $T$ is finite hence $\sharp{\cP} < \infty$ and the union on the right of~\eqref{eq_preimage} is finite. The set $\tilde{K}$ is a compact, so it is closed.

We remark that we incidentally verified that the space $\tilde{D}$ is compact.

\subsection{}\label{vlast_kompakt_pr} The space $D$ is a compactum.

$D$ is the compact space as a factor-space of compact space $\tilde{D}$. Compactum $\tilde{D}$ is the normal topological space and a factor-space of a normal space over closed partition is a normal space, see~\cite{RF}. Thus $D$ is a normal space, in particularly, $D$ is Hausdorff space. Therefore $D$ is compactum.

\subsection{}\label{vlast_vkladennya_T} Map $\varphi = \restrict{\pr}{T} : T \rightarrow D$ is the embedding.

By definition, $F_{x} \cap T = \{x\}$ for every $x \in T$, hence $\varphi$ is an injective map. The space $T$ is compact and the space $D$ is Hausdorff thus $\varphi$ is homeomorphism onto its image, see~\cite{RF}.

\subsection{}\label{vlast_vkladennya_diskiv} For every $P \in \cP$ a map $\restrict{\pr}{D_{P}} : D_{P} \rightarrow D$ is an embedding.

By definition, for $x \in D_{P}$ we get
\begin{equation}
D_{P} \cap F_{x} =
\left\{
\begin{array}{ll}
	\Phi_{P}(\Phi_{P}^{-1}(x)) \,, & x \in \Phi_{P}(P) \,, \\
	\{ x \} \,, & x \in D_{P} \setminus \Phi_{P}(P) \,.
\end{array}
\right.
\end{equation}\label{eq_intersect_disk}
  The map $\Phi_{P}$ is injective hence $\Phi_{P}(\Phi_{P}^{-1}(x)) = \{ x \}$, $x \in \Phi_{P}(P)$. Finally, $F_{x} \cap D_{P} = \{ x \}$ for every $x \in D_{P}$ and a continuous map $\restrict{\pr}{D_{P}}$ is injective. Thus it is a homeomorphism of compact $D_{P}$ onto its image.

\subsection{}\label{vlast_vnutrishnist_obraza_diska} For every $P \in \cP$ a set $\pr(D_{P} \setminus \Phi_{P}(P))$ is open in $D$ and has no common points with a set $\pr(\tilde{D} \setminus (D_{P} \setminus \Phi_{P}(P)))$.

 Let $P \in \cP$. The set $D_{P}$ is open-closed in the space $\tilde{D}$, hence an open set  $D_{P} \setminus \Phi_{P}(P)$ in  $D_{P}$  is also open in $\tilde{D}$. This set is saturated by definition. Therefore $D_{P} \setminus \Phi_{P}(P) = \pr^{-1}(\pr(D_{P} \setminus \Phi_{P}(P)))$ and a set $\pr(D_{P} \setminus \Phi_{P}(P))$ is open in the factor-space $D$.

It follows from the discussion above that a set $\tilde{D} \setminus (D_{P} \setminus \Phi_{P}(P))$ is also saturated and it has no common points with $D_{P} \setminus \Phi_{P}(P)$. Thus
\[
\pr(D_{P} \setminus \Phi_{P}(P)) \cap \pr(\tilde{D} \setminus (D_{P} \setminus \Phi_{P}(P))) = \emptyset \,.
\]

\subsection{}\label{vlast_okil_rebra_disk} Let $e \in E$ be any edge of the tree $T$, points $w_{1}$, $w_{2}$ be the ends of $e$ and $P'$, $P'' \in \cP$ be paths in $\cP$ that pass through $e$. We designate $e^{0} = e \setminus \{ w_{1}, w_{2} \}$,
\begin{eqnarray*}
D^{0}_{P'} & = & D_{P'} \setminus \partial D_{P'} \subseteq \bigcup_{p \in \cP} D_{P} \setminus \Phi_{P}(P) \,, \\
D^{0}_{P''} & = & D_{P''} \setminus \partial D_{P''} \subseteq \bigcup_{p \in \cP} D_{P} \setminus \Phi_{P}(P) \,, \\
\tilde{U} & = & (D^{0}_{P'} \cup \Phi_{P'}(e^{0})) \sqcup (D^{0}_{P''} \cup \Phi_{P''}(e^{0})) \sqcup e^{0} \,, \\
U & = & \pr(\tilde{U}) \,.
\end{eqnarray*}
 $U$ is the open neighborhood of a set $\pr(e^{0})$ in the space $D$, it is homeomorphic to open disk and is divided by a set $\pr(e^{0})$ onto two connected components $\pr(D^{0}_{P'})$ and $\pr(D^{0}_{P''})$.

To prove this we should remark that sets $e^{0}$, $D^{0}_{P'}$ and $D^{0}_{P''}$ are open in $\tilde{D}$. By definition of partition $\ff$ for every $x \in e^{0}$ we get $F_{x} = \{ x, \Phi_{P'}(x), \Phi_{P''}(x) \}$ since the set $\tilde{U}$ is saturated.

The set $\tilde{U}$ is open in $\tilde{D}$. Really, in the first place $e^{0}$ is an open subset of $T$, secondly, $\Phi_{P'}(e^{0})$ is an open subset of closed subspace $\Phi_{P'}(P')$ of space $D_{P'}$, therefore, $\Phi_{P'}(P') \setminus \Phi_{P'}(e^{0})$ is a closed subset $D_{P'}$. Let us remark that arcs $\Phi_{P'}(P')$ and $\partial D_{P'} \setminus \Phi_{P'}(P')$ have $\Phi_{P'}$- images of endpoints of the path $P'$as common ends, thus $\Cl{\partial D_{P'} \setminus \Phi_{P'}(P')} \cap \Phi_{P'}(e^{0}) = \emptyset$ and following conditions hold true
\begin{gather*}
\partial D_{P'} \setminus \Phi_{P'}(e^{0}) = (\partial D_{P'} \setminus \Phi_{P'}(P')) \cup (\Phi_{P'}(P'))  \setminus \Phi_{P'} (e^{0})) = \\
= \Cl{(\partial D_{P'} \setminus \Phi_{P'}(P'))} \cup (\Phi_{P'}(P'))  \setminus \Phi_{P'}(e^{0})) \,.
\end{gather*}
So, a set $\partial D_{P'} \setminus \Phi_{P'}(e^{0})$ is closed in $D_{P'}$ and a set $D^{0}_{P'} \cup \Phi_{P'}(e^{0}) = D_{P'} \setminus (\partial D_{P'} \setminus \Phi_{P'}(e^{0}))$ is open in $D_{P'}$. Similarly, a set $D^{0}_{P''} \cup \Phi_{P''}(e^{0})$ is open in $D_{P''}$. Sets $T$, $D_{P'}$ and $D_{P''}$ are open-closed in space $\tilde{D}$. Thus the set $\tilde{U}$ is open in $\tilde{D}$.

Finally, the set $U = \pr(\tilde{U})$ is open in $D$. This set is a result of gluing
\begin{gather*}
U \cong (D^{0}_{P''} \cup \Phi_{P''}(e^{0})) \cup_{\alpha} (D^{0}_{P'} \cup \Phi_{P'}(e^{0})) \,, \\
\alpha = \Phi_{P''} \circ \Phi_{P'}^{-1} : \Phi_{P'}(e^{0}) \rightarrow \Phi_{P''}(e^{0}) \,.
\end{gather*}

 A map $\alpha$ is a composition of homeomorphisms. Therefore $U$ is homeomorphic to open disk  and is divided by $\pr(e^{0})$ onto two connected components $\pr(D^{0}_{P'})$ and $\pr(D^{0}_{P''})$.

\subsection{}\label{vlast_mezsha} For any $P_{1}, \ldots, P_{n} \in \cP$ a boundary $\Fr{D_{n}}$ of a set $D_{n} = \pr(\bigcup_{i=1}^{n} D_{P_{i}})$ in the space $D$ belongs to $\pr(\bigcup_{i=1}^{n} P_{i}) = \pr(T) \cap D_{n}$.

It follows from property~\ref{vlast_vnutrishnist_obraza_diska} that $\Fr{\pr(D_{P_{i}})} \subseteq \pr(\Phi_{P_{i}}(P_{i})) = \pr(P_{i})$ for any $i \in \{1, \ldots, n\}$. Hence
\[
\Fr{D_{n}} \subseteq \bigcup_{i=1}^{n} \Fr{\pr(D_{P_{i}})} \subseteq \bigcup_{i=1}^{n} \pr(P_{i}) = \pr\Bigl(\bigcup_{i=1}^{n} P_{i}\Bigr) \,.
\]

\subsection{}\label{vlast_rebro_i_diski} Let $P_{1}, \ldots, P_{n} \in \cP$, $\tilde{D}_{n} = \bigcup_{i=1}^{n} D_{P_{i}}$, $D_{n} = \pr(\tilde{D}_{n})$. Let $e$ be an edge of $T$ such that $\pr(e^{0}) \cap D_{n} \neq \emptyset$, where $e^{0}$ is an edge $e$ without ends.

 A set $\pr(e^{0})$ belongs to $\Int{D_{n}}$ iff at least one point $y \in \pr(e^{0})$ has a neighborhood in $D_{n}$ which is homeomorphic to open disk. Otherwise, a set  $\pr(e^{0})$ belongs to $\Fr{D_{n}}$.

 If $\pr(e^{0}) \subseteq \Int{D_{n}}$ then a set $\pr(e^{0})$ has a neighborhood in $D_{n}$ which is homeomorphic to open disk and both paths $P'$, $P'' \in \cP$ passing through $e$ belong to a set $\{P_{1}, \ldots, P_{n}\}$.

If $\pr(e^{0}) \subseteq \Fr{D_{n}}$, then exactly one of them belongs to $\{P_{1}, \ldots, P_{n}\}$.


Suppose that paths $P'$, $P'' \in \cP$ pass through the edge $e$. By the definition $\pr(e^{0}) \subseteq D_{n} \cap \pr(T) = \pr(\bigcup_{i=1}^{n} P_{i})$, so  at least one of them belongs to $\{P_{1}, \ldots, P_{n}\}$. We consider two possibilities.

We assume that $P' = P_{k}$, $P'' = P_{s}$, $k, s \in \{1, \ldots, n\}$. Then a set
\[
U = \pr((D^{0}_{P'} \cup \Phi_{P'}(e^{0})) \cup (D^{0}_{P''} \cup \Phi_{P''}(e^{0})) \cup e^{0}) \subseteq D_{n}
\]
is an open neighborhood of $\pr(e^{0})$ that is homeomorphic to an open disk, see~\ref{vlast_okil_rebra_disk}.

Let $P' \in \{P_{1}, \ldots, P_{n}\}$, $P'' \notin \{P_{1}, \ldots, P_{n}\}$. In this case $U = U' \cup U'' \cup e^{0}$, $U' = \pr(D^{0}_{P'}) \subseteq D_{n}$ but a set $U'' = \pr(D^{0}_{P''})$ has no common points with $D_{n}$, therefore $\pr(e^{0}) \subseteq \Fr{D_{n}}$.

Suppose that for some $y \in \pr(e^{0})$ in $D_{n}$ there exists an neighborhood $W_{y} \in D_{n}$ homeomorphic to open two dimensional disk. By using theorem of Shenflies~\cite{Newman, Z-F-C} we can find a small neighborhood $\hat{W}_{y}$ of $y$ in $D_{n}$ such that it is homeomorphic to an open disk and satisfies following conditions:
\begin{itemize}
	\item A set $\Cl{\hat{W}_{y}}$ in the space $D_{n}$ is homeomorphic to a closed disk and is separated from compacts $\pr(T \setminus e^{0})$ and $D \setminus U$.
	\item $\Cl{\hat{W}_{y}}$ intersects $\pr(e^{0})$ by a connected segment that is a cut of the disk $\Cl{\hat{W}_{y}}$.
\end{itemize}
Then the set $\pr(e^{0})$ divides $\hat{W}_{y}$ onto two connected components $W_{1} \cup W_{2} = \hat{W}_{y} \setminus \pr(e^{0})$, $W_{1} \cap W_{2} = \emptyset$ such that $\Cl{W}_{1} \cap \Cl{W}_{2} \ni y$.

 By the construction $\hat{W}_{y} \subseteq U \cap D_{n}$ and $W_{1}$, $W_{2} \subseteq (U \cap D_{n}) \setminus \pr(T)$. Let us remind that $P' \in \{P_{1}, \ldots, P_{n}\}$, thus $\pr(D_{P'}) = D' \subseteq D_{n}$. Similarly, $P'' \notin \{P_{1}, \ldots, P_{n}\}$ hence $\pr(D^{0}_{P''}) \cap D_{n} = \emptyset$, see property~\ref{vlast_vnutrishnist_obraza_diska}. Therefore $U \cap D_{n} = U' \cup \pr(e^{0})$, where $U' = \pr(D^{0}_{P'})$, see property~\ref{vlast_okil_rebra_disk}, and $U \cap D_{n} \cap \pr(T) = \pr(e^{0}) \subseteq \partial D'$, where $D' = \pr(D_{P'})$.

Thus $y \in \partial D'$ and the set $\hat{W}_{y}$ is the open neighborhood of $y$ in closed disk $D'$ and $\hat{W}_{y} \cap (D' \setminus \partial D') = W_{1} \cup W_{2}$, $y \in \Cl{W}_{1} \cap \Cl{W}_{2}$. By Lemma~\ref{lem_loc_connect} we can conclude that $W_{1} = W_{2}$ but it contradicts to the assumption that $W_{1} \cap W_{2} = \emptyset$.

So, if $\{P', P''\} \nsubseteq \{P_{1}, \ldots, P_{n}\}$, then there is no $y \in \pr(e^{0})$ that has an open neighborhood in $D_{n}$, which is homeomorphic to open disk.

\subsection{}\label{vlast_struct_mezshi} Let $P_{1}, \ldots, P_{n} \in \cP$. Let us describe a structure of boundary $\Fr{D_{n}}$ of $D_{n} = \pr(\bigcup_{i=1}^{n} D_{P_{i}})$ in $D$.

 Denote by $E_{n} \subseteq E$ a set of all edges of the tree $T$ such that exactly one of two paths $P'$, $P'' \in \cP$ passing through $e \in E_{n}$ belongs to $\{P_{1}, \ldots, P_{n}\}$. As we know, see Condition~\ref{vlast_rebro_i_diski}, $\pr(E_{n}) \subseteq \Fr{D_{n}}$ and if for some edge $e \in E$ we get $e \notin E_{n}$, then $\Fr{D_{n}} \cap \pr(e) \subseteq \{v', v''\}$, where $v'$, $v'' \in V$ are ends of $e$.

Similarly, denote by $V_{n} \subseteq V$ a set of all vertices of $T$ such that for a vertex $v \in V_{n}$ the following condition satisfies: $\pr(v) \in D_{n}$ and all edges that are adjacent to $v$ belong to $E \setminus E_{n}$. It is easy to show that the set $V_{n}$ is discreet and $\pr(E_{n}) \cap \pr(V_{n}) = \emptyset$.

From the discussion above and Condition~\ref{vlast_mezsha} it follows that
\begin{equation}\label{eq_struct_mezshi}
\pr(E_{n}) \subseteq \Fr(D_{n}) \subseteq (\pr(E_{n}) \cup \pr(V_{n})) \,.
\end{equation}

\subsection{}\label{vlast_zvyazny_pidgraf} Let $P_{1}, \ldots, P_{n} \in \cP$. A set $D_{n} = \pr(\bigcup_{i=1}^{n} D_{P_{i}})$ is connected iff then $\bigcup^{n}_{i=1} P_{i}$ is a connected subgraph of the tree $T$.

 Let $\bigcup^{n}_{i=1} P_{i} = T'$ is a connected subgraph of $T$. Then $D_{n} = \pr(T') \cup \bigcup_{i=1}^{n} \pr(D_{P_{i}})$, all sets $\pr(T')$, $\pr(D_{P_{i}})$, $i \in \{1, \ldots, n\}$ are connected and $\pr(T') \cap \pr(D_{P_{i}}) \neq \emptyset$, $i \in \{1, \ldots, n\}$. Hence the set $D_{n}$ is connected.

Next, let $\bigcup_{i=1}^{n} D_{P_{i}} = T' \cup T''$, $T' \cap T'' = \emptyset$ and sets $T'$, $T''$ are nonempty and closed. Every set $P_{i}$, $i \in \{1, \ldots, n\}$ is connected, therefore, either $P_{i} \in T'$ or $P_{i} \in T''$. Without loss of generality we can change indexing of the elements of  $\{P_{1}, \ldots, P_{n}\}$ in such way that for some $s \in \{1, \ldots, n-1\}$ the following conditions are satisfied
\[
T' = \bigcup_{i=1}^{s} P_{i} \,, \quad T'' = \bigcup_{i=s+1}^{n} P_{i} \,.
\]

Every set
\[
\tilde{D}' = T' \cup \bigcup_{i=1}^{s} D_{P_{i}} \,, \quad \tilde{D}'' = T''\bigcup_{i=s+1}^{n} D_{P_{i}} \,,
\]
is closed, whence sets $D' = \pr(\tilde{D}')$ і $D'' = \pr(\tilde{D}'')$ are closed, see Condition~\ref{vlast_zamkn_pr}. By the construction $\tilde{D}' \cap \tilde{D}'' = \emptyset$. Let $y \in D' \cap D''$. A map $\pr$ is injective by definition on the set $\pr^{-1}(D \setminus \pr(T))$ and sets $\tilde{D}'$ and $\tilde{D}''$  do not intersect on $\pr^{-1}(D \setminus \pr(T))$, thus $y \in \pr(T)$. Hence $y \in \pr(T \cap \tilde{D}') \cap \pr(T \cap \tilde{D}'') = \pr(T') \cap \pr(T'')$. But as we know, see Condition~\ref{vlast_vkladennya_T}, the map $\varphi = \restrict{\pr}{T}$ is bijective, therefore $\pr(T') \cap \pr(T'') = \pr(T' \cap T'') = \emptyset$. We get a contradiction, thus  $D' \cap D'' = \emptyset$.

Hence $D_{n} = D' \sqcup D''$ an sets $D'$, $D''$ are closed and nonempty. Therefore the set $D_{n}$ is not connected.

\medskip

Finally let us prove a $\cD$-planarity of the tree $T$.	

Let for some $n$, $1 \leq n < \sharp \cP$ directed paths $P_{1} = P(v_{1}, v_{1}'), \ldots, P_{n} = P(v_{n}, v_{n}') \in \cP$ are fixed and $\tilde{D}_{n} = \bigcup_{i=1}^{n} D_{P_{i}}$, $D_{n} = \pr(\tilde{D}_{n})$.

For every $i \in \{1, \ldots, n\}$ we denote by $\tilde{\gamma}_{i}$ an directed arc of $\partial D_{P_{i}}$ from point $\Phi_{P_{i}}(v_{i}')$ to $\Phi_{P_{i}}(v_{i})$ which has no other common points with an arc $\Phi_{P_{i}}(P_{i})$.

Suppose that the objects under consideration comply with following conditions.
\begin{itemize}
	\item [(i)] A space $D_{n}$ is homeomorphic to a close two-dimensional disk.
	\item [(ii)] There exists at least one edge $e \in \bigcup_{i=1}^{n} P_{i}$ such that its image $\pr(e)$ is contained in a boundary circle $\partial D_{n}$ of $D_{n}$.
	\item [(iii)] A disk $D_{n}$ is oriented in the following way: for every $i \in \{1, \ldots, n\}$ and every edge $e \in P_{i}$ such that $\pr(e)$ belongs to $\partial D_{n}$ an orientation of $e$ generated by the direction of $P_{i} = P(v_{i}, v_{i}')$ maps by $\pr$ onto an orientation of $D_{n}$.
	\item [(iv)] For every $i \in \{1, \ldots, n\}$ an arc $\gamma_{i} = \pr(\tilde{\gamma}_{i})$ connects a point $\pr(v_{i}')$ with a point $\pr(v_{i})$ and has no other common points with a set $\pr(T)$ and orientation of this arc is consistent with the orientation of $D_{n}$.
\end{itemize}

We should remark that for $n = 1$ and any path $P = P_{1} \in \cP$ if we take an orientation on $D_{1}=\pr(D_p)$ induced from $D_{P}$ by using $\pr$, then  Conditions (i)--(iv) always hold true. By the construction, Conditions (iii) and (iv) are true , (i) follows from Condition ~\ref{vlast_vkladennya_diskiv}, (ii)  follows from Condition~\ref{vlast_rebro_i_diski}.

We also remark that it follows from Condition~\ref{vlast_rebro_i_diski} that an edge $e \in \bigcup_{i=1}^{n} P_{i}$ belongs to $\partial D_{n}$ of $D_{n}$ iff $e \in E_{n}$. Thus Condition~(iii) is well-posed. As well all boundary points of $D_{n}$ in the space $D$ possibly except a finite number of isolated points from the set $\pr(V_{n})$ belong to $\partial D_{n}$.

Let an edge $e \in \bigcup_{i=1}^{n} P_{i}$ satisfies Condition~(ii). Then $e \in E_{n}$ and there is the unique path $P_{n+1} = P(v_{n+1}, v_{n+1}') \in \cP \setminus \{P_{1}, \ldots, P_{n}\}$ such that it passes through the edge $e$. Let $e \in P_{l}$, where $P_{l} \in \{P_{1}, \ldots, P_{n}\}$ is the second path among two paths from the set $\cP$ which passes through the edge $e$.

Let us consider a disk $D_{P_{n+1}}$ and its image $D' = \pr(D_{P_{n+1}})$. By Condition~\ref{vlast_vkladennya_diskiv} it is also the closed disk. Let $\Gamma = D_{n} \cap D'$. It is obvious that $\Gamma$ is closed.

By Condition~\ref{vlast_vnutrishnist_obraza_diska} a set $\pr(D_{P_{n+1}} \setminus \Phi_{P_{n+1}}(P_{n+1}))$ is open in $D$ and does not intersect $D_{n}$. It follows from Condition~\ref{vlast_vkladennya_diskiv} that
\begin{gather*}
\pr(D_{P_{n+1}} \setminus \Phi_{P_{n+1}}(P_{n+1})) = \pr(D_{P_{n+1}}) \setminus \pr \circ \Phi_{P_{n+1}}(P_{n+1}) = D' \setminus \pr(P_{n+1}) \,, \\
\pr(\Phi_{P_{n+1}}(P_{n+1})) = \pr(P_{n+1}) \subseteq \Cl{(D' \setminus \pr(P_{n+1}))} \,.
\end{gather*}
Therefore
\[
\Gamma = \Fr{D_{n}} \cap \Fr{D'} \subseteq \pr(P_{n+1}) \,.
\]

Let us apply Condition~\ref{vlast_struct_mezshi} to $D_{n}$ and $D'$. By~\eqref{eq_struct_mezshi} the set $\Gamma$ consists of images of edges which belong to the path $P_{n+1}$ and possibly from a number of images of vertices of a tree $T$.

Let us check that the set $\Gamma$ is connected.

If it is not the case it follows from what we said above that there are two vertices $w_{1}$, $w_{2} \in V$, $w_{1} \neq w_{2}$ of $T$ such that they belong to the path $P_{n+1}$ and a projection of a path $P(w_{1}, w_{2}) \subseteq P_{n+1}$ which connects them in $T$ intersects $\Gamma$ by a set $\{ \pr(w_{1}), \pr(w_{2}) \}$. Then $\pr(P(w_{1}, w_{2})) \cap D_{n} = \{w_{1}, w_{2}\}$.

On the other hand, the set $D_{n}$ is connected thus $T' = \bigcup_{i=1}^{n} P_{i}$ is a connected subgraph of $T$, see Condition~\ref{vlast_zvyazny_pidgraf}. From Condition~\ref{vlast_mezsha} it follows that $w_{1}$, $w_{2} \in V(T')$, therefore there is a path $P'(w_{1}, w_{2})$ connecting them in $T'$. This path has to connect $w_{1}$ with $w_{2}$ in $T$. But $\pr(P'(w_{1}, w_{2})) \subseteq D_{n}$ hence $P'(w_{1}, w_{2}) \neq P(w_{1}, w_{2})$. So, vertices $w_{1}$ and $w_{2}$ of $T$ can be connected in $T$ by two different paths which is impossible in the tree $T$.

This contradiction proves that $\Gamma$ is connected.

It follows from the connectedness of $\Gamma$ and from the inclusion $\pr(e) \subseteq \Gamma \cap \pr(E_{n})$ that $\Gamma \subseteq \pr(E_{n})$. Thus
\[
\Gamma \subseteq \partial D_{n} \cap \partial D' \,,
\]
where $\partial D' = \pr(\partial D_{P_{n+1}})$ is a boundary circle of the disk $D'$.

By the discussion above and from $\Gamma \subseteq  \pr(P_{n+1})$ it is easy to understand that
\[
\Gamma = \pr(P(v, v'))
\]
for some $v$, $v' \in V \cap P_{n+1}$, $v \neq v'$.

It is obvious that $P(v, v')$ is homeomorphic to a closed segment. From the Conditions~\ref{vlast_vkladennya_T} and~\ref{vlast_vkladennya_diskiv} it follows that it is embedded into a boundary circles $\partial D_{n}$ and $\partial D'$ by means of maps
\begin{gather*}
\psi_{n} = \restrict{\pr}{P(v, v')} : P(v, v') \rightarrow D_{n} \,, \\
\psi' = \pr \circ \Phi_{P_{n+1}} : P(v, v') \rightarrow D' \,.
\end{gather*}
Therefore, a set
\[
D_{n+1} = D_{n} \cup D' \cong D_{n} \cup_{\psi} D' \,, \quad \psi = \psi_{n} \circ (\psi')^{-1} \,,
\]
is a result of a gluing of closed disks $D_{n}$ and $D'$ by a segment that is embedded into the boundary circles of these disks. Consequently the set $D_{n+1}$ is homeomorphic to a closed disk.

Let us denote $\tilde{D}_{n+1} = \bigcup_{i=1}^{n+1} D_{P_{i}}$. It is clear that
\[
D_{n+1} = \pr\Bigl(\bigcup_{i=1}^{n} D_{P_{i}}\Bigr) \cup \pr(D_{P_{n+1}}) = \pr\Bigl(\bigcup_{i=1}^{n+1} D_{P_{i}}\Bigr) = \pr(\tilde{D}_{n+1}) \,.
\]
Hence the space $D_{n+1}$ constructed according to the set $\{P_{1}, \ldots, P_{n+1}\}$ satisfies Condition~(i).

Disks $D_{n}$ and $D'$ are oriented. The orientation of $D'$ is generated by an orientation of $D_{P_{n+1}}$ by means of the map $\pr$.

By Condition~(iii) applied to $D_{n}$ and $D'$ we get two orientations on $e$. One of them is induced from an orientation of $P_{l} \supseteq e$ and is coordinated with orientation of $D_{n}$. Another is generated by direction of $P_{n+1}$ and is consistent with an orientation of $D'$.

As we said above the directed paths $P_{l}$, $P_{n+1} \in \cP$ containing an edge $e$ have to pass through $e$ in the opposite directions. Therefore the orientations induced on $\Gamma$ from $D_{n}$ and $D'$ are opposite. Hence the orientations of $D_{n}$ and $D'$ are coordinated and generate an orientation of $D_{n+1}$. It complies with the following condition
\begin{itemize}
	\item for any simple arc $\alpha : I \rightarrow \partial D_{n} \cap \partial D_{n+1}$ an orientation of $\alpha$ is consistent with orientation of $D_{n+1}$ iff  an orientation $\alpha$ is coordinated with orientation of $D_{n}$;
	\item for any simple arc $\beta : I \rightarrow \partial D' \cap \partial D_{n+1}$ an orientation of $\beta$ is consistent with an orientation of $D_{n+1}$ iff it is coordinated with an orientation a disk $D'$.
\end{itemize}

Disks $D_{n}$ and $D'$ satisfy Conditions~(iii) and~(iv). So, according to what has being said  $D_{n+1}$ also satisfies Conditions~(iii) and~(iv).

Suppose that the set $D_{n+1}$ does not satisfy Condition~(ii). Then $E_{n+1} = \emptyset$, see Condition~\ref{vlast_struct_mezshi} and Remark ~(iii), and $\partial D_{n+1} \cap \pr(T) \subseteq \pr(V)$. Thus a set $\partial D_{n+1} \cap \pr(T)$ is finite.

The following correlations are implicated from  Condition~\ref{vlast_vnutrishnist_obraza_diska}
\[
\partial D_{n+1} \setminus \pr(T) \subseteq \pr\Bigl( \bigcup_{i=1}^{n+1} (D_{P_{i}} \setminus \Phi_{P_{i}}(P_{i})) \Bigr) = \bigcup_{i=1}^{n+1} \pr(D_{P_{i}} \setminus \Phi_{P_{i}}(P_{i})) \,.
\]
From Condition~\ref{vlast_vkladennya_diskiv} it follows that for every $i \in \{1, \ldots, n+1\}$ a set $\pr(D_{P_{i}} \setminus \partial D_{P_{i}}) \subseteq D_{n+1}$ is homeomorphic to an open disk. Hence
\[
\bigcup_{i=1}^{n+1} \pr(D_{P_{i}} \setminus \partial D_{P_{i}}) \subseteq D_{n+1} \setminus \partial D_{n+1} \,.
\]
From this correlation it follows, see Condition~(iv), that
\begin{gather*}
\partial D_{n+1} \setminus \pr(T) \subseteq \left[ \bigcup_{i=1}^{n+1} \bigl(\pr(D_{P_{i}} \setminus \partial D_{P_{i}}) \cup \pr(\partial D_{P_{i}} \setminus \Phi_{P_{i}}(P_{i}))\bigr) \right] \cap \partial D_{n+1} = \\
= \bigcup_{i=1}^{n+1} \pr(\partial D_{P_{i}} \setminus \Phi_{P_{i}}(P_{i})) \subseteq \bigcup_{i=1}^{n+1} \pr(\tilde{\gamma}_{i}) = \bigcup_{i=1}^{n+1} \gamma_{i} \,.
\end{gather*}

 A set $\bigcup_{i=1}^{n+1} \gamma_{i}$ is closed in $D$ hence it is also closed in $\partial D_{n+1}$. Therefore, a set $\partial D_{n+1} \setminus \bigcup_{i=1}^{n+1} \gamma_{i}$ have to be an open subset of a space $\partial D_{n+1}$. But
\[
\partial D_{n+1} \setminus \bigcup_{i=1}^{n+1} \gamma_{i} \subseteq \partial D_{n+1} \cap \pr(T) \subseteq \pr(V)
\]
and this set is finite. Consequently,
\[
\partial D_{n+1} = \bigcup_{i=1}^{n+1} \gamma_{i} \,.
\]

 From Condition~(iv) it easily follows that open arcs $\gamma_{i} \setminus \{\pr(v_{i}), \pr(v_{i}')\}$, $i \in \{1, \ldots, n+1\}$ are pairwise disjoint. Therefore every point of a set $\partial D_{n+1} \cap \pr(T) = \bigcup_{i=1}^{n+1} \{\pr(v_{i}), \pr(v_{i}')\}$ is a common boundary point of exactly two arcs of the family $\{\gamma_{i}\}_{i=1}^{n+1}$.

It follows from the choice of an orientation of arcs $\gamma_{i}$, $i \in \{1, \ldots, n+1\}$ that if for some $s$, $r \in \{1, \ldots, n+1\}$ either $v_{s} = v_{r}$ or $v_{s}' = v_{r}'$ is true, then $s=r$. Thus for every $i \in \{1, \ldots, n+1\}$ there is the unique  $j(i) \in \{1, \ldots, n+1\}$, such that $v_{i} = v_{j}'$ and if $r \neq s$ then $j(r) \neq j(s)$ . We also remark that by the construction $n \geq 1$, thus $n+1 \geq 2$ and $j(i) \neq i$, $i \in \{1, \ldots, n+1\}$.

Therefore, on the set $\{1, \ldots, n+1\}$ there is a transposition $\sigma$ without fix points such that $v_{i} = v_{\sigma(i)}'$, $i \in \{1, \ldots, n+1\}$. Let $\sigma = c_{1} \cdots c_{k}$ be a decomposition of $\sigma$ into independent cycles. Let $c_{1} = (i_{1} \ldots i_{m})$. Then $v_{i_{1}} = v_{i_{2}}', \ldots v_{i_{m-1}} = v_{i_{m}}'$, $v_{i_{m}} = v_{i_{1}}'$.

From the definition of the set $\cP$ we get $v_{i}' \rho_{C} v_{i}$, $i \in \{1, \ldots, n+1\}$, since $P_{i} = P(v_{i}, v_{i}') \in \cP$. So, it is true that
\[
v_{i_{1}} \rho_{C} v_{i_{2}},\; \ldots,\; v_{i_{m-1}} \rho_{C} v_{i_{m}},\; v_{i_{m}} \rho_{C} v_{i_{1}} \,,
\]
thus vertices of the set $M_{1} = \{v_{i_{1}}, \ldots, v_{i_{m}}\}$ generate a $\rho_{C}$-cycle, see Definition~\ref{ozn_cycle}. From Corollary~\ref{nas_struct_2} it is follows that the set $M_{1}$ is a class of equivalence of a minimal equivalence relation $\hat{\rho}_{C}$ which contains the relation $\rho_{C}$. By Proposition~\ref{prop_cycle} and Corollary ~\ref{nas_struct_2} the relation $\hat{\rho}_{C}$ has the unique class of equivalence $V^{\ast}$. Hence $M_{1} = V^{\ast}$, $\sigma = c_{1}$, $n+1 = \sharp V^{\ast} = \sharp \cP$ and $D_{n+1} = D$.

From what was said above it follows that for $n+1 < \sharp \cP$ the disk $D_{n+1}$ satisfies Condition~(ii). Thus, for $n+1 < \sharp \cP$ the disk $D_{n+1}$ satisfies~(i)--(iv), but for $n+1 = \sharp \cP$ it complies with conditions~(i) and~(iv).

Finally, starting from any path $P = P_{1} = P(v_{1}, v_{1}') \in \cP$, we can sort out elements of a set
\[
\cP = \{P_{1}=P(v_{1}, v_{1}'), \ldots, P_{N}=P(v_{N}, v_{N}')\}
\]
in a finite number of steps so that for every set
\[
D_{n} = \pr\Bigl( \bigcup_{i=1}^{n} D_{P_{i}} \Bigr) \,, \quad n \in \{1, \ldots, N-1\} \,,
\]
the conditions~(i)--(iv) are true and for the set
\[
D_{N} = \pr\Bigl( \bigcup_{i=1}^{N} D_{P_{i}} \Bigr) = \pr\Bigl( \bigcup_{P \in \cP} D_{P} \Bigr) = \pr(\tilde{D}) = D
\]
conditions~(i) і~(iv) are also true.

Thus $D_{N} = D$ is closed oriented two-dimensional disk, $\varphi = \restrict{\pr}{T} : T \rightarrow D$ is an embedding, see Condition~\ref{vlast_vkladennya_T}.

 For every edge $e \in E$ both paths of $\cP$ passing through this edge belong to a set $\{P_{1}, \ldots, P_{N}\}$, thus $E_{N} = \emptyset$ and $\partial D = \bigcup_{i=1}^{N} \gamma_{i}$ with open arcs $\gamma_{i} \setminus \{\pr(v_{i}), \pr(v_{i}')\}$ are pairwise disjoint. It is clear that
\[
\varphi(T) \cap \partial D = \bigcup_{i=1}^{N} \; \{\pr(v_{i}), \pr(v_{i}')\} = \bigcup_{P(v, v') \in \cP} \{\pr(v), \pr(v')\} = V^{\ast} \,.
\]

An orientation of $D$ generates some cyclic order $O$ on the set $\pr(V^{\ast})$ . A map $\varphi_{0} = \restrict{\varphi}{V^{\ast}} : V^{\ast} \rightarrow \pr(V^{\ast})$ is bijective, therefore, a map $\varphi_{0}^{-1}$ generates on the set $V^{\ast}$ some cyclic order $C'$ which is an isomorphic image of a cyclic order $O$ ($C'(v_{1}, v_{2}, v_{3}) \Leftrightarrow O(\pr(v_{1}), \pr(v_{2}), \pr(v_{3}))$).

We induce a convenient relation  $\rho_{C'}$ on $V^{\ast}$, see Definition~\ref{ozn_induc_by_cycle}. From Condition~(iv) it follows that for every $i \in \{1, \ldots, N\}$ we have $v_{i}' \rho_{C'} v_{i}$. On the other hand, by definition of the set $\cP$ it follows that $v' \rho_{C} v$ iff $P(v, v') \in \cP$. But $\cP = \{P_{1}, \ldots, P_{N}\}$, hence if $P(v, v') \in \cP$, then $P(v, v') = P_{i} = P(v_{i}, v_{i}')$ for some $i \in \{1, \ldots, N\}$. Therefore the following conditions hold true
\[
v' \rho_{C} v \Rightarrow v' \rho_{C'} v \,, \quad v, v' \in V^{\ast} \,,
\]
and the relation $\rho_{C'}$ contains $\rho_{C}$.

With the help of convenient relations $\rho_{C}$ and $\rho_{C'}$ we can induce on $V^{\ast}$ the relations of cyclic orders $C_{\rho_{C}}$ and $C_{\rho_{C'}}$, respectively, see Definition~\ref{ozn_induc_by_zruchn} and Proposition~\ref{prop_induced_cyclic}. From Definition~\ref{ozn_induc_by_zruchn} it is easily  follows that if $\rho_{C'}$ contains $\rho_{C}$ then $C_{\rho_{C'}}$ contains $C_{\rho_{C}}$. In other words, an identical map $Id_{V^{\ast}}$ is monomorphism of cyclic order $C_{\rho_{C}}$ onto $C_{\rho_{C'}}$, see Definition~\ref{ozn_cyclic_morphisms}. From  Lemma~\ref{lem_induced_induced_equiv} it follows that $C_{\rho_{C}} = C$ and $C_{\rho_{C'}} = C'$, hence the map  $Id_{V^{\ast}}$ is monomorphism of the cyclic order $C$ onto $C'$. Lemma~\ref{lem_full_cycle_iso} implies that the map $Id_{V^{\ast}}$ is an isomorphism of cyclic order $C$ onto $C'$.

By the construction a map $\varphi_{0}^{-1}$ is an isomorphism of cyclic order $O$ onto $C'$ thus $\varphi_{0}$ is an isomorphism of cyclic order $C = C'$ onto a cyclic order $O$ which is induced onto $\varphi(V^{\ast})$ from an oriented circle $\partial D$.

Finally, the map $\varphi$ satisfies all conditions of Definition ~\ref{ozn_01} and a tree $T$ is $\cD$-planar.

\end{proof}

 

\end{document}